\documentclass[12pt]{article} 
\usepackage{amsmath,amssymb,amsthm,eucal,MnSymbol,amscd,mathbbol,tikz} 
\usepackage{hyperref}
\usepackage{float}
\usepackage{multirow}
\usepackage{xfrac}
\usepackage{pifont}
\usepackage{subcaption}

\font\tenrm=cmr10

\font\cmssl=cmss10 at 12 pt  
   
\font\bigss=cmssdc10 scaled 2300

\font\cmsslll=cmss10 at 14 pt  
\newcommand{\g}{\gamma}  
  
\newcommand\Sp{\mathrm{Sp}}

\newcommand\SU{\mathrm{SU}}  
\newcommand\U{\mathrm{U}}

\newcommand{\ra}{\rightarrow}

\DeclareMathOperator\Inn{Inn}

\DeclareMathOperator\Out{Out}
\DeclareMathOperator\Aut{Aut}
\DeclareMathOperator\Real{Re}
\DeclareMathOperator\Imag{Im}
\DeclareMathOperator\Ker{Ker}
\DeclareMathOperator\BT{2T}
\DeclareMathOperator\BO{2O}
\DeclareMathOperator\BI{2I}
\DeclareMathOperator\BD{2D}
\DeclareMathOperator\D{D}
\DeclareMathOperator\T{T}
\DeclareMathOperator\OR{O}
\DeclareMathOperator\I{I}
\DeclareMathOperator\Alt{Alt}
\DeclareMathOperator\Sym{Sym}

\DeclareMathOperator\ID{Id}

\newtheorem{Th}{Theorem}[section]  
\newtheorem{Prop}[Th]{Proposition}  
\newtheorem{Cor}[Th]{Corollary}  
\newtheorem{Lem}[Th]{Lemma}  

\theoremstyle{definition} 
\newtheorem{Def}[Th]{Definition}  
\newtheorem{Ex}[Th]{Example}

\theoremstyle{remark}

\newtheorem*{remark}{Remark}
\newtheorem*{cav}{Caveat}

\newcommand{\bt}{\begin{Th}\ \ }  
\newcommand{\et}{\end{Th}}  
\newcommand{\bp}{\begin{Prop}\ \ }  
\newcommand{\ep}{\end{Prop}}  
\newcommand{\bc}{\begin{Cor}\ \ }  
\newcommand{\ec}{\end{Cor}}  
\newcommand{\bl}{\begin{Lem}\ \ }  
\newcommand{\el}{\end{Lem}}  
\newcommand{\bd}{\begin{Def}\ \ }  
\newcommand{\ed}{\end{Def}}  
\newcommand{\bex}{\begin{Ex}\ \ }  
\newcommand{\eex}{\end{Ex}}  

\newcommand{\pf}{\begin{proof}}
\newcommand{\epf}{\end{proof}}

\newcommand{\be}{\begin{equation}}  
\newcommand{\ee}{\end{equation}}

\newcommand{\arr}{\begin{array}{rlll}}  
\newcommand{\ea}{\end{array}}  
\newcommand{\bea}{\begin{eqnarray}}  
\newcommand{\eea}{\end{eqnarray}}  
\newcommand{\bean}{\begin{eqnarray*}}  
\newcommand{\eean}{\end{eqnarray*}}  
  
\catcode`@=11  
\@addtoreset{equation}{section}  
\catcode`@=12

\begin{document}  

\vskip 1.5 true cm  
\begin{center}  
{\bigss   Locally homogeneous nearly K\"ahler manifolds}  

\vspace{2ex}

\vskip 1.0 true cm   
{\cmsslll  V.\ Cort\'es$^1$ and J.\ J.\ V\'asquez$^2$} \\[3pt] 
{\tenrm   $^1$Department Mathematik
und Zentrum f\"ur Mathematische Physik\\ 
Universit\"at Hamburg, 
Bundesstra{\ss}e 55, 
D-20146 Hamburg, Germany\\  
cortes@math.uni-hamburg.de}\\[1em]   
{\tenrm   $^2$Max-Planck-Institut f\"ur Mathematik
in den Naturwissenschaften\\ 
Inselstra{\ss}e 22, D-04103 Leipzig, Germany\\
Jose.Vasquez@mis.mpg.de}\\[1em] 
October 24, 2014 
\end{center}  
\vskip 1.0 true cm  
\baselineskip=18pt  
\begin{abstract}  
We construct locally homogeneous  6-dimensional nearly K\"ahler manifolds as 
quotients of homogeneous nearly K\"ahler manifolds $M$ by freely acting finite subgroups of $\Aut_0(M)$.
We show that non-trivial such groups do only exists if $M=S^3\times S^3$. In that case we classify all 
freely acting subgroups of $\Aut_0(M)=\SU (2) \times \SU (2) \times  \SU (2)$ of the form $A\times B$, where $A\subset
\SU (2) \times \SU (2)$ and $B\subset \SU (2)$.  
\noindent  
\end{abstract}
\tableofcontents

\section*{Introduction}
Recall that an almost Hermitian manifold $(M,g,J)$ without nontrivial K\"ahler local de Rham factor is called  {\cmssl (strict) nearly K\"ahler} if 
$(\nabla_X J)X=0$ for all $X\in TM$, where $\nabla$ denotes the Levi-Civita connection.   
It was shown by Nagy \cite{N} that all complete simply connected nearly K\"ahler manifolds 
are products of twistor spaces of quaternionic K\"ahler manifolds of positive scalar curvature, 
homogeneous spaces and six-dimensional nearly K\"ahler manifolds. 
 
According to Butruille \cite{B1,B2} there exist only 4 examples
of $6$-dimensional homogeneous nearly K\"ahler manifolds $M=G/K$: 
\begin{enumerate}
\item the sphere $S^6=\mathrm{G}_2/{\SU(3)}$,
\item the complex projective space $\mathbb{C}P^3 = \Sp (2)/(\U (1) \times \Sp(1))$, 
\item the flag manifold ${F}_{1,2}(\mathbb{C}^3)=\SU(3)/(\U(1)\times U(1))$,
\item the Lie group  $S^3\times S^3 = \SU(2)^3/\Delta(\SU (2))$, where $\Delta : \SU(2) 
\hookrightarrow \SU(2)^3$ is the diagonal embedding.
\end{enumerate} 
To our knowledge, these exhaust all examples of $6$-dimensional nearly K\"ahler manifolds which have 
occurred in the literature so far. Incidentally, the second and third examples are precisely the twistor spaces
of the $4$-dimensional quaternionic K\"ahler manifolds of positive scalar curvature. 
Each of these four homogeneous spaces $M=G/K$ is a $3$-symmetric space and $G=\mathrm{Aut}_0 (M)$ is the maximal connected group of automorphisms of the 
nearly K\"ahler structure. The latter statement follows from \cite[Theorem 5.3]{GM}, which uses \cite[Theorem 3.6]{T}. 

In this paper we are interested in six-dimensional nearly K\"ahler manifolds $M$ for which the 
universal covering $\tilde{M}$ is homogeneous.  Such manifolds will be called 
{\cmssl locally homogeneous nearly K\"ahler manifolds} in the following. 
The classification of these manifolds amounts to the description 
of the finite subgroups $\Gamma \subset \mathrm{Aut}(\tilde{M})$ which act 
freely on $\tilde{M}$, for each of the $4$ (simply connected) homogeneous nearly K\"ahler manifolds $\tilde{M}$ from
Butruille's list.  For simplicity, we will only consider subgroups $\Gamma$ of $G=\mathrm{Aut}_0 (\tilde{M})$. 
The corresponding locally homogeneous nearly K\"ahler manifolds $M$ 
are the quotients $M=\tilde{M}/\Gamma = \Gamma \backslash G/K$, by the natural left-action of $\Gamma\subset G$ on $\tilde{M}=G/K$. 

The next proposition shows that it is sufficient to consider the case $\tilde{M} = S^3 \times S^3$, in which case 
we will classify certain classes of freely acting groups of automorphisms in the main part of the paper.
\bp Let $M$ be a homogeneous nearly K\"ahler manifold such that $G = \mathrm{Aut}_0(M)$ admits 
a nontrivial subgroup acting freely on $M$. Then the nearly K\"ahler manifold $M$ is isomorphic to $S^3 \times S^3$. 
\ep
\pf 
Any element $\gamma \in G$ is contained in some maximal torus $T$ of $G$. If 
$M$ is not isomorphic to $S^3 \times S^3$ then the stabilizer $K\subset G$ of a point $o\in M$ is of maximal rank and, hence, contains a maximal 
torus $T_0$ of $G$. Since any two maximal tori are conjugate, there exists an element 
$a\in G$ such that $aTa^{-1} = T_0$. This implies that $p=a^{-1}o\in M$ is a fixed point of $\g$.  This shows that 
$G$ does not contain any nontrivial subgroup $\Gamma$ acting freely on $M$. 
\epf
From now on we consider the case $M=S^3\times S^3=G/K$, where  $G=\SU (2) \times \SU (2) \times  \SU (2)$ and 
$K=\Delta(\SU (2))\subset G$. Notice that the nearly K\"ahler structure on $M$ can be considered as a left-invariant 
structure on the Lie group $L=\SU (2) \times \SU (2)=M$. Let $L$ act by left-translations as a 
subgroup of $G=\mathrm{Aut}_0 (M)$, where the inclusion is simply $(a,b) \mapsto (a,b, \mathbb{1})$. Since the action of $L=\SU (2) \times \SU (2) \cong \SU (2) \times \SU (2) \times \{ \mathbb{1} \} \subset G$ by left-translations on $M$ is free, any finite subgroup $\Gamma \subset L$ gives rise to a locally homogeneous nearly K\"ahler manifold $M/\Gamma$. 

 Our main results amount to the classification of all finite subgroups of $G=\Aut_0(M)=L\times \SU (2)$ acting freely on $M$ that (up to a permutation of the three factors of $G$)  are of the form $A\times B$, for some finite subgroups $A\subset L$, $B\subset \SU (2)$. In addition, we classify all finite simple groups $\Gamma \subset G$ that act freely on $M$, see Theorem \ref{simple} in Section 2.  We refer to Theorem \ref{main} in Section 3 for the description of finite subgroups $\Gamma_1 \times \Gamma_2 \times \Gamma_3 \subset G $ acting freely on $M$ that are products of groups $\Gamma_i \subset \SU(2)$ for $i=1,2,3$, and to Theorems  \ref{typeB}, \ref{type3}, \ref{q.final}, \ref{q.final2} in Section 4, for the remaining groups. This yields a wealth of new examples of  nearly K\"ahler manifolds.
 
 {\em Acknowledgements:} This work was supported by the Collaborative Research Center SFB 676 ``Particles, Strings, and the Early Universe'' of the Deutsche Forschungsgemeinschaft. 

\section{Goursat's Theorem}
Finite subgroups of the product $G_1 \times G_2$ of two abstract groups are described by Goursat's theorem, see e.g. \cite{FIG}.  
We give the proof of the theorem as it is needed in the sequel. 
\begin{Th}[Goursat's theorem] \label{GOURSAT}
Let $G_1, G_2$ be groups. There is a one-to-one correspondence between subgroups $C  \subset G_1 \times G_2$ and quintuples $\mathcal{Q}(C)=\{A,A_0,B,B_0,\theta\}$, where $A_0 \triangleleft A \subset G_1$, $B_0 \triangleleft B  \subset G_2 $ and $\theta: \sfrac{A}{A_0} \longrightarrow \sfrac{B}{B_0} $ is an isomorphism.
\end{Th}
\pf Let $C \subset G_1\times G_2$ be a subgroup and denote by $\pi_i: G_1\times G_2 \longrightarrow G_i$, $i=1,2$,  the natural projections.  Set $A= \pi_1 (C) \subset G_1$, $B=\pi_2(C)\subset G_2$, $A_0= \Ker(\pi_2|_{_{C}})$ and $B_0=\Ker(\pi_1|_{_{C}})$. It is readily seen that $A_0$ and $B_0$ can be identified with normal subgroups of $A$ and $B$ respectively. We denote these subgroups again by $A_0$ and $B_0$. Define a map $\tilde{\theta}:A \longrightarrow \sfrac{B}{B_0}$ as follows. For $a\in A$ pick any $b \in B$ so that $(a,b)\in C$ and set $\tilde{\theta}(a):=bB_0$. One can check that this map is well-defined and factorizes through an isomorphism $\theta: \sfrac{A}{A_0} \longrightarrow \sfrac{B}{B_0}$. This defines a map $C \mapsto \mathcal{Q}(C)$.

Conversely, a quintuple $Q=\{A,A_0,B,B_0, \theta \}$ as described above defines a group $C=\mathcal{G}(Q)  \subset G_1 \times G_2$ by setting $C=p^{-1}(\Gamma(\theta))$, where $p: A\times B \longrightarrow \sfrac{A}{A_0} \times \sfrac{B}{B_0}$  is the natural  homomorphism and $\Gamma( \theta) \subset \sfrac{A}{A_0} \times \sfrac{B}{B_0}$ denotes the graph of the homomorphism $\theta$. Observe that $C \subset G_1 \times G_2$ is in fact a fiber product, 
\begin{align}
 C& = \{ (a,b)\in A \times B \ : \ \theta (aA_0)=bB_0\} = \{ (a,b)\in A\times B| \; \alpha (a) = \beta (b)\} , \label{fiber-prod.1} 
\end{align}
where  
\begin{center}
\begin{tikzpicture}[node distance=2cm, auto]
  \node (A) {$\alpha:A$};
	\node (B) [right of=A] {$\sfrac{A}{A_0}$};
  \node (C) [right of=B]{$\sfrac{B}{B_0}$};  
  \node (D) [right of=C]{ and};
   \node (E) [right of=D]{$\beta:B $};
    \node (F) [right of=E]{$ \sfrac{B}{B_0}$};
  \draw[->] (A) to node {} (B);
  \draw[->] (E) to node {} (F);
	\draw[->] (B) to node {$\theta$} (C);
	  \end{tikzpicture}
\end{center}
are the natural homomorphisms. The maps $\mathcal{Q}$ and $\mathcal{G}$ are inverse to each other.
 \epf
\bp
Two subgroups $C,C' \subset G_1 \times G_2$  with corresponding quintuples $\mathcal{Q}(C)=\{A,A_0,B,B_0,\theta\}$, 
$\mathcal{Q}(C')=\{A',A'_0,B',B'_0,\theta'\}$
are conjugate if and only if  there exists $(g_1,g_2)\in G_1 \times G_2$ such that 
$A'= g_1 A g^{-1}_1, B'= g_2 B g^{-1}_2, A'_0= g_1 A_0 g^{-1}_1, B'_0= g_2 B_0 g^{-1}_2$ and the diagram 
\begin{center}
\begin{tikzpicture}[node distance=2cm, auto]
  \node (A) {$A/A_0$};
\node (B) [right of=A] {$B/B_0$};
\node (C) [below of= A] {$A'/A'_0$};
\node (D) [below of= B] {$B'/B'_0$};
---------
\draw[->] (A) to node [above] {$\theta $} (B);
\draw[->] (A) to node [swap] {$c(g_1)$} (C);
\draw[->] (B) to node [right]   {$c(g_2)$} (D);
\draw[->] (C) to node [swap] {$\theta'$} (D);
  \end{tikzpicture}
  \end{center}  
  commutes. Where $c(g_i)$ denotes  conjugation by $g_i \in G_i$, $i=1,2$. 
\ep

\begin{remark} \label{Remark} 
Sometimes we will consider different subgroups $C=\mathcal{G}(A,A_0,B,B_0,$ $\theta )\subset G_1\times G_2$
for fixed $A,A_0,B,B_0$.  In that case it is convenient to identify $A/A_0$ and $B/B_0$ with the same 
abstract group $F$ and consider $\theta : A/A_0 \ra B/B_0$ as an automorphism of $F$. 
\end{remark} 

In the remaining parts of this paper, we classify (up to conjugation) finite subgroups  $C \subset G$  acting freely on the nearly K\"ahler manifold $M$, which are either simple or of the form $ D \times E \subset  G$, for $D \subset \SU(2)^2$ and $E \subset \SU(2)$ arbitrary. This motivates the following definition.
\begin{Def}
A   finite subgroup $C \subset G= \SU(2)^3 $ is said to be  \emph{splittable} whenever $C=A_1 \times A_2 \times A_3 \subset G$ for some non-trivial  subgroups  $A_i \subset \SU(2)$ for $i=1,2,3$, and  \emph{semi-splittable} if $C= D \times E$ for some non-trivial subgroups $D \subset \SU(2)^2$, $E \subset \SU(2)$. In addition, a semi-splittable group $C \subset G$ is said to be  \emph{strict} if it is not splittable. 
\end{Def} 
 We are excluding the case of trivial factors in the above definition because the occurrence of a trivial factor 
 implies that $C$ acts freely on $M$, as mentioned in the introduction.

\section{Simple groups}
The Lemma below will help us to distinguish subgroups of $G=\SU (2) \times \SU (2) \times \SU (2)$ that act freely on $M=S^3\times S^3$. 
\begin{Lem}\label{crit}
A subgroup $C \subset G$ acts non freely on $M$ if and only if there is a non-trivial element $(a_1, a_2,a_3) \in C $, so that $$ \Real (a_1)= \Real (a_2)  =\Real (a_3).$$
\end{Lem}
\pf Consider the action of $\SU(2)=S^3\subset \mathbb{H} = \mathbb{R}^4$ on itself given by conjugation. The orbit of a given unit quaternion $1\neq a=\Real(a)+\Imag(a)$ is then of the form 
\begin{align}
\{ \Real (a)\} \times S^2(\rho) \subset \mathbb{R}\times \mathbb{R}^3= \mathbb{R}^4  & \ \ , \  \ \ \rho = \sqrt{1-\Real (a)^2}. \label{cd}  
\end{align}
Let $(a_1,a_2,a_3)\in C$ be a non-trivial element fixing a class  $[(a,b,c)]\in M$, i.e.\  so that $(a,b,c)^{-1} (a_1,a_2,a_3)(a,b,c)=(w,w,w) \in K$. As the real part of a quaternion is invariant under conjugation, we have $\Real (a_1)=\Real(a_2)=\Real(a_3)=\Real(w)$. Conversely, suppose there is an element $(1,1,1)\neq (a_1,a_2,a_3)$ $ \in C$ with: $\Real(a_1)=\Real(a_2)=\Real(a_3)$. Then, by relation (\ref{cd}) the orbits of $a_1,a_2,a_3 \in \SU(2)$ are the same.\epf
Let us recall now that any finite subgroup of $\SU(2)$ is conjugate to one of the so-called ADE groups, see e.g. Theorem 1.2.4 in \cite{GAB}. These groups are described in terms of generators as follows.  
\begin{table}[H]
\begin{center}
  \begin{tabular}{ ||c | c |c | c || }
    \hline
		Label & Name &  Order & Generators \\ \hline
			$\mathbb{A}_{n-1}$ & $\mathbb{Z}_n$ & $n$ & $e^{\frac{i 2 \pi}{n} }$ \\ 
			$\mathbb{D}_{n+2}$ & $\BD_{2n}$ & $4n$ & $j, e^{\frac{i \pi }{n}}$ \\ 
			$\mathbb{E}_6$ & $\BT$ & 24 & $\frac{1}{2}(1+i)(1+j), \frac{1}{2}(1+j)(1+i)$ \\ 
			$\mathbb{E}_7$ & $\BO$ & 48 & $\frac{1}{2}(1+i)(1+j), \frac{1}{\sqrt{2}}(1+i)$ \\ 
			$\mathbb{E}_8$ & $\BI$ & 120 & $\frac{1}{2}(1+i)(1+j), \frac{1}{2}(\phi+\phi^{-1}i +j)$ \\ \hline
\end{tabular}
\end{center}
\caption{Finite subgroups of $\SU(2)$.}
\label{tabla1}
\end{table}
Here $n \geq 2$ and $\phi= \frac{1+\sqrt{5}}{2}$ is the golden ratio. An element-wise description of the ADE groups is provided in Table \ref{tablax}.
\begin{table}[H]
\begin{center}
  \begin{tabular}{ ||c | c |c | c || }
    \hline
		Label & Name &  Order & Elements \\ \hline
			$\mathbb{A}_{n-1}$ & $\mathbb{Z}_n$ & $n$ & $\{e^{\frac{ 2 \pi i x}{n} } \ : \ x=0,...,n-1 \}$ \\ 
			$\mathbb{D}_{n+2}$ & $\BD_{2n}$ & $4n$ & $\{e^{\frac{i  \pi x}{n} } :  x=0,...,2n-1 \} \cup j \{e^{\frac{i  \pi x}{n} }  :  x=0,...,2n-1 \} $ \\ 
			$\mathbb{E}_6$ & 2T & 24 & $\BD_{4}\cup \left\{ \frac{\pm 1 \pm i \pm j \pm k}{2}\right\}$ \\ 
			$\mathbb{E}_7$ & 2O & 48 & $\BT \cup e^{\frac{i\pi}{4}}\BT$ \\ 
			$\mathbb{E}_8$ & 2I & 120 & $\BT \cup q \BT \cup q^2 \BT \cup q^3 \BT \cup q^4 \BT$ \\ \hline
\end{tabular}
\end{center}
\caption{Element description for ADE groups.}
\label{tablax}
\end{table}
Where  $q= \frac{1}{2}\left(\phi +\phi^{-1}i +j \right)$.  
\begin{Th}\label{simple}
The following are, up to conjugation and permutation of the factors, the only non-trivial simple subgroups $C \subset G$ acting freely on $M$. 
\begin{enumerate}
\item[(a)] $\mathbb{Z}_p \times \{1\} \times \{ 1 \}$.
\item[(b)] $\Gamma(\varphi(r)) \times \{ 1\}$, where $\varphi(r) \in \Aut(\mathbb{Z}_p) $, see Section \ref{automorphisms}. 
\item[(c)] $C(p,r,s)=\{(x,\varphi(r)x,\varphi(s)x) \ : x \in \mathbb{Z}_p \}$, where $\varphi(r),\varphi(s)$ are automorphisms of $\mathbb{Z}_p$ and either $r \neq \pm 1 \mod p$, or $s \neq \pm 1 \mod p$. 
\end{enumerate} 
Here $p\in \mathbb{Z}$ is an arbitrary prime number. 
\end{Th}
\pf Let $C \subset G$ be a non trivial simple subgroup and $j\in \{1,2,3\}$  be so that $\pi_j(\Gamma)\neq \{1\}$, where $\pi_j$ denotes the natural projection. Because $C$ is simple, the taken projection restricts to an isomorphism $\left(\pi_j\right)|_{C}: C  \longrightarrow \pi_j(C) \subset G$. However, groups of type  DE are non commutative and have non trivial center, whilst groups of type A are commutative and non simple, unless they are of prime order. Consequently, $C \subset G$ has the isomorphism type of $\mathbb{Z}_p$ for some fixed prime number $p$. We distinguish between the following cases 
\begin{itemize}
	\item[(a)] Let $\pi_1(C) \neq \{1\}$ and $\pi_2(C)=\pi_3(C)=\{1\}$. That is, $C= A \times \{1\}\times \{1\}$ for some finite subgroup $A \subset \SU(2)$. Because $\pi_1 |_{C}$ an isomorphism, the group $C \subset G$ is conjugate to $\mathbb{Z}_p \times \{1\}\times \{1\} $. 
	\item[(b)] Let $\pi_1(C), \  \pi_2(C)\neq \{1\}$ and $\pi_3(C)= \{1\}$. In particular  $C=D \times \{1\} \subset G$ for some finite subgroup $D \subset \SU(2)^2$. Denote by $\mathcal{Q}(D)=\{A, A_0, B, B_0, \theta\}$ the quintuple defining $D\subset \SU(2)^2$, see Theorem \ref{GOURSAT}. By construction of $\mathcal{Q}(D)$ and the  simplicity of $C$, we see that $A_0=B_0=\{1\}$ and (up to conjugation) $A=B= \mathbb{Z}_p \subset \SU(2)$, and so, the isomorphism $\theta$ can be realized as an automorphism of $\mathbb{Z}_p$. It follows that $C\subset  G$ is conjugate to $\Gamma(\theta) \times \{1\}$.
		\item[(c) ] Let $\pi_j(C)\neq \{1\}$ for $j=1,2,3$, $G_1= \SU(2)^2$ and $G_2= \SU(2)$. The subgroup $C \subset  G_1 \times G_2$ determines the quintuple $\mathcal{Q}(C)=\{A,A_0,B,B_0, \theta\}$, where 
		\begin{align}
		A=&\{ (a,b) \in G_1 \ : \ (a,b,c)\in C \ \textrm{for some } c\in \SU(2) \} , \nonumber \\
        B=& \{ c \in G_2 \ : \  (a,b,c) \in C \ \textrm{for some } c\in \SU(2) \} , \nonumber \\		
        A_0=&\{ (a,b) \in A \ : \ (a,b,1)\in C \} , \nonumber \\		
		B_0=&\{ c\in B \ : \ (1,1,c)\in C \} , \nonumber 
		\end{align}
		and $\theta: \sfrac{A}{A_0}\rightarrow \sfrac{B}{B_0}$ is an isomorphism. Using the simplicity of $C \subset G$, we can easily check that $A$ is simple, $A_0=\{(1,1)\}$ and $B_0=\{1\}$. Just as in (b), we conclude that (possibly after conjugation) $A =\Gamma(\varphi) \subset G_1$ for some $\varphi\in \Aut(\mathbb{Z}_p)$. The isomorphism $\theta$ can be thus realized as an automorphism of $\mathbb{Z}_p$, and so, the group $C \subset  G$ is conjugate to a group of the form $\{ (x,\varphi (x),\psi (x)) :  x \in \mathbb{Z}_p\}$, where $\varphi, \psi$ are automorphism of $\mathbb{Z}_p$.		
\end{itemize}
We are now left to decide using Lemma \ref{crit}, which of the groups described above act freely on $M$. To this end, observe that the groups $\mathbb{Z}_p \times \{1\} \times \{ 1 \}$ and $\Gamma(\varphi(r)) \times \{ 1\}$ act freely on $M$, so we can suppose $C=C(p,r,s) \subset G$ for integers $r,s \in \mathbb{Z}$ so that $r,s \nmid  p$. Lemma \ref{crit} tells us that the group $C(p,r,s)\subset  G$ acts freely on $M$ precisely if the system
\begin{align}
(1\pm r)x  = 0 \mod p & \ , \ rx = \pm sx \mod \ p \label{system}  
\end{align}
admits just the trivial solution. 
\epf
\section{Splittable groups} 
In this section we will classify freely acting splittable subgroups $A_1 \times A_2 \times A_3 \subset  G= \SU(2)^3$. The correspondence given in Theorem \ref{GOURSAT} together with our knowledge of ADE groups will be used in order to construct all relevant subgroups $C \subset  \SU(2)^2$. Thereafter, to verify using Lemma \ref{crit} when a given subgroup $A_1 \times A_2 \times A_3\subset G$ acts in the desired fashion amounts to solve certain kinds of equations on integers, such as in (\ref{system}), to which we devote the forthcoming section.
\subsection{Integral equations}
To begin with, we summarize all real parts of elements in ADE groups.          
\begin{table}[H] 
\begin{center}
  \begin{tabular}{ ||c | c |c || }
    \hline
		Label 			 & Name 	&   Real parts \\ \hline
			$\mathbb{A}_{n-1}$ & $\mathbb{Z}_n$ &   $\cos\left(\frac{2 \pi x}{n}\right) \ \ 0 <x \leq n$ \\ 
			$\mathbb{D}_{n+2}$ & $\BD_{2n}$ 			& 	$0 ,\cos\left(\frac{\pi x}{n}\right) \ \ 0 < x \leq 2n $ \\ 
			$\mathbb{E}_6$ 		 & $\BT$  			 			&	  $0, \pm 1 , \pm \frac{1}{2}$ \\ 
			$\mathbb{E}_7$ 		 & $\BO$ 			 			&   $0 , \pm 1 , \pm \frac{1}{2} , \pm \frac{\sqrt{2}}{2} , \pm \frac{1}{2\sqrt{2}}$  \\ 
			$\mathbb{E}_8$								& $\BI$ & $0, \pm 1 , \pm \frac{1}{2} ,  \frac{\pm 1 \pm \sqrt{5} }{4}$ \\ \hline
\end{tabular}
\end{center}
\caption{Real parts of elements in ADE groups.}
\label{tabla4}
\end{table}
We are interested in solving the following type of equations on integers:

\underline{Case I} $ax+by=c$, for $a,b \in \mathbb{Z}  \setminus  \{0\}, c \in \mathbb{Z}$. This case corresponds to a linear Diophantine equation which can be completely solved using Bezout's  lemma, see e.g. Theorems 2.1.1 and 2.1.2 in \cite{DIO}.
\begin{Th} \label{lind}
Let $a,b \in \mathbb{Z} \setminus  \{0\}$ and $c \in \mathbb{Z}$. The equation 
\begin{align}
ax+by&=c \label{lin}
\end{align}
is solvable if and only if $\gcd(a,b)| c$, in which case its general solution reads 
\begin{align}
(x,y)&=(x_0,y_0)+ t(b,-a)  \ \ \ t \in \mathbb{Z}.
\end{align} 
Here $(x_0,y_0)\in \mathbb{Z}^2$ is a particular solution of equation \eqref{lin}.
\end{Th}
\underline{Case II} $\cos\left(\frac{2\pi x}{n}\right)= b, \  b\in \left\{ 0, \pm 1 , \pm \frac{1}{2}, \pm \frac{1}{\sqrt{m}} \ : \ m\in \mathbb{N}\right\}$. We determine first all possible values of $m \in \mathbb{N}$  so that $b=\pm \frac{1}{\sqrt{m}}$ can be written as a cosine of some rational multiple of $2 \pi$, by making use of the following standard result in algebraic number theory, see e.g. Section 3 of \cite{COS}. 
 \begin{Th}
Let $\theta$ be a rational multiple of $2 \pi$. If $\cos(\theta)\in \mathbb{Q}$, then $\cos(\theta) \in \left\{ -1, - \frac{1}{2}, 0 , \frac{1}{2}, 1 \right\}$.
\end{Th}
In our situation, it follows that $m\in \{1,2,4\}$, and so, $b \in \mathbb{R}$ belongs to the following list of values: $0, \pm 1 , \pm \frac{1}{\sqrt{2}}, \pm \frac{1}{2}$. Table \ref{tabla2} displays the restrictions on $n \in \mathbb{Z}$ in order to get an integral solution $x(n)\in \mathbb{Z}$ for the equation $\cos\left(\frac{2\pi x}{n}\right)= b$.
\begin{table}[H]
\begin{center}
  \begin{tabular}{ ||c | c |c || }
    \hline
		$b$					 &   $x(n)\in \mathbb{Q}$ & Restriction on $n$ \\ \hline
				$-1$									 &			$\frac{n}{2}+nk$																	&  $2|n$							\\ 
		$- \frac{1}{\sqrt{2}}$ &	   $\frac{3n}{8}+nk, \frac{5}{8}n+nk$ &$8 |n$\\  
					$- \frac{1}{2}$  &     $\frac{n}{3}+nk \ , \ \frac{2n}{3}+nk$ & $3 | n$  \\ 
					$0$ 						 &	  $\frac{n}{4}+nk \ , \ \frac{3}{4}n+nk$ & $4 | n$ \\ 
					$\frac{1}{2}$		 &    $\frac{n}{6}+nk \ , \ \frac{5n}{6}+nk$ &$6 |n$\\   
					$\frac{1}{\sqrt{2}}$&  $\frac{n}{8}+nk \ , \ \frac{7n}{8}+nk$ & $8 | n $\\  
					$1$							&		$nk$ & -			 \\ \hline
\end{tabular}
\end{center}
\caption{Conditions on $n$.}
\label{tabla2}
\end{table}
\underline{Case III} $\cos\left(\frac{2\pi x}{n}\right)=c, \ c\in \{ a+\sqrt{5}b$ \ : \  $a,b\in \mathbb{Q}^{\times} \}$. Let $c=a+b \sqrt{5}$ with $a,b\in \mathbb{Q}^{\times}$ and suppose $x(n)\in \mathbb{Z}$ solves the equation in consideration. From this setup, we can easily derive that $x(n)\in \mathbb{Z} $ satisfies the integral equation $\cos(2\theta) +B \cos(\theta) = C$, for $B=-4a$ and $C=10b^2-1-2a^2$, and where $\theta=\frac{2 \pi x}{n}$. The latter equation can be completely solved by means of Theorem 7 in \cite{CON}, result which is stated in generality sufficient to our needs. 
\begin{Th}[Conway, Jones]
Suppose we have at most two distinct rational multiples of $\pi$ lying strictly between $0$ and $\frac{\pi}{2}$ for which some rational linear combination of their cosines is rational. Then the appropriate linear combination is proportional to either $\cos(\sfrac{\pi}{5})-\cos(\sfrac{2\pi}{5})=\frac{1}{2}$ or $\cos(\sfrac{\pi}{3})=\frac{1}{2}$. 
\end{Th}
Table \ref{tabla3} displays the additional values of $c\in \mathbb{Q}[\sqrt{5}]$ we must consider, together with the restrictions to impose on $n\in \mathbb{Z}$ in order to find an integral solution $x(n) \in \mathbb{Z}$ for the equation $\cos\left(\frac{2\pi x}{n}\right)=c$.  
\begin{table}[H]
\begin{center}
  \begin{tabular}{ ||c | c  | c|| }
    \hline
		$c= a +b \sqrt{5}$					 &    $x(n) \in \mathbb{Q}$  & Restriction on $n$		\\ \hline
		      $\frac{1}{4}(1+\sqrt{5})$ &	    $\frac{n}{10}+kn, \frac{9n}{10}+kn$ &	$10 | n$									\\ 
					$\frac{1}{4}(-1+\sqrt{5})$ &     $\frac{n}{5}+kn, \frac{4n}{5}+kn$ &			$5|n$							\\ 
					$\frac{1}{4}(1-\sqrt{5})$ &	    $\frac{3n}{10}+kn, \frac{7n}{10}+kn$ &			$10|n$						\\ 
					
					$-\frac{1}{4}(1+\sqrt{5})$ &     $\frac{2n}{5}+kn, \frac{3n}{5}+kn$& 		$5|n$								\\ \hline  
					 
\end{tabular}
\end{center}
\caption{Conditions on $n$.}
\label{tabla3}
\end{table}
\begin{remark}
In the latter situation $\cos(\theta)$ is an algebraic number of degree 2. It is in fact a zero of the quadratic polynomial $(t-a)^2-5b^2$. The classification of algebraic numbers of low degree and sufficiently small length, see e.g. Section 5 of \cite{COS}, provides an alternative argument to build Table \ref{tabla3}. 
\end{remark}
\subsection{Automorphisms of quotient groups}\label{automorphisms} 
The present section comprises descriptions of automorphisms groups of quotients of ADE groups that are relevant in the forthcoming sections. The main reference is Section 6.2 in \cite{FIG}. This material will be however adapted to our needs.

(1) The group of outer automorphisms of $\mathbb{Z}_n$ is given by 
\begin{align}
\Out(\mathbb{Z}_n)&=\{ \varphi(r) \ : \ \gcd (r,n)=1 \} , \nonumber
\end{align}
where $\varphi(r)$ denotes the map $\mathbb{Z}_n \ni x \mapsto x^r \in \mathbb{Z}_n$ in multiplicative notation. 

(2) To describe the outer automorphism group of a dihedral group $\D_{2n}$, consider the following presentation of $\D_{2n}$  
\begin{align}
\D_{2n}&=\left< x,y \ : \ x^2=y^n=(xy)^2=1\right> = \{ y^p  : \ 0 \leq p < n\} \cup \{ xy^p  : \ 0 \leq p < n\} . \nonumber
\end{align}
 Observe $\D_2 = \mathbb{Z}_2$, so we can assume that $n>1$. The case $n=2$ is also special as $\D_4$ is isomorphic to the Klein  Vierergruppe. The automorphism group of $\D_4$ isomorphic to $\Sym(3)$ and acts by permutations of the 3 non trivial involutions. The outer automorphism group of $\D_{2n}$ for $n >2$  is 
		\begin{align}
		\Out (\D_{2n})&=\left< \tau_{a,b} \ : \ (a,b)\in \mathbb{Z}^{\times}_{n} \times \mathbb{Z}_{n} \right> \cong \mathbb{Z}^{\times}_{n} \ltimes \mathbb{Z}_{n},  \nonumber
		\end{align}
	where the action of the affine group $\mathbb{Z}^{\times}_{n} \ltimes \mathbb{Z}_{n}$ on $\D_{2n}$ is given by
	\begin{align}
	\tau_{a,b}(y^p)&=y^{ap} \ \ , \ \ \tau_{a,b}(xy^p)= xy^{ap+b}. \nonumber  
	\end{align}
(Here $\mathbb{Z}_{n}$ denotes the additive group and $\mathbb{Z}_{n}^{\times}$ the multiplicative group of units in the ring $\mathbb{Z}_{n}$.) 
	
(3) Since $\BD_2 \cong \mathbb{Z}_4$,  we consider $\BD_{2n}$ only for $n>1$. 
We have the following  presentation: 
\begin{align}
\BD_{2n}&=\left< s,t  : \ s^2=t^n=(st)^2 \right> =\{ t^p   : \ 0 \leq p <2n \} \cup \{ st^p  \ : \ 0 \leq p <2n \} . \nonumber
\end{align} 
In fact, we can take $s=je^{i\frac{\pi}{n}}$ and $t=e^{i\frac{\pi}{n}}$ when $\BD_{2n}$ is realized as a subgroup
of $\SU (2)$. 
The outer automorphism group of $\BD_{2n}$ for $n>2$ is also an affine  group:
\begin{align}
		\Out (\BD_{2n})&=\left< \tau_{a,b} \ : \ (a,b)\in \mathbb{Z}^{\times}_{2n} \times \mathbb{Z}_{2n} \right> 
		 \cong \mathbb{Z}^{\times}_{2n} \ltimes \mathbb{Z}_{2n}, \nonumber  
		\end{align}
	where the action on $\BD_{2n}$ is given by
	\begin{align}
	\tau_{a,b}(t^p)&= t^{ap} \ \ , \ \	\tau_{a,b}(st^p)= st^{ap+b}. \nonumber
	\end{align}	
	We need to make a distinction for $n=2$. Any automorphism of $\BD_4=\{\pm 1, \pm i , \pm j ,  \pm k\} \subset \SU(2)$ is obtained via conjugation with an element in $\BO$ modulo $\mathbb{Z}_2=\{\pm1\}$. The point-wise action of $\sfrac{\BO}{\mathbb{Z}_2}$ on $\BD_4$ is described below.

\begin{table}[H]
\begin{center}
 \begin{tabular}{ c || c | c | c | c |c || c |c |c}   
						&  $i$ & $j$ &  $k$  & & &  $i$ & $j$ &  $k$ \\ \hline \hline	
			$[i]$		&   $i$ &  $-j$ & $-k$ &      &  $[\frac{1}{\sqrt{2}}(1-i)]$	 & $i $ & $-k $&  $j $ 	\\
		    $[j]$		 & $-i$ & $j$ & $-k$ &        & $[\frac{1}{\sqrt{2}}(j+k)]$	   & $-i $  & $k $&  $j $ 	\\
			$[k]$		 & $-i$  & $-j$ & $k$ 	&     	& $[\frac{1}{\sqrt{2}}(j-k)]$	 & $-i $ & $-k $ & $-j $\\
			 $[\frac{1}{2}(1+i+j+k) ]$  & $ j$  & $ k$ & $ i$ & & $[\frac{1}{\sqrt{2}}(i+k)]$	   & $ k$  & $-j $&  $i $ 	\\
			$[\frac{1}{2}(1-i-j-k)]$   & $ k$  & $ i$& $ j$ & &      $[\frac{1}{\sqrt{2}}(1-k)]$	  & $ -j$  & $i $&  $k $ 	\\

		    $[\frac{1}{2}(1+i-j-k)]$   & $ -j$ & $ k$&  $ -i$	& \quad \quad & 		$[\frac{1}{\sqrt{2}}(i-k)]$	   & $-k $  & $-j $&  $-i $ 	\\
		    $[\frac{1}{2}(1+i+j-k)]$ & $ -k$  & $ i$&  $ -j$ 	& & 		$[\frac{1}{\sqrt{2}}(i+j)]$	  & $j $  & $i $&  $-k $ 	\\	
		    $[\frac{1}{2}(1-i+j-k)]$ & $ -j$  & $ -k$& $ i$ 	& & 		$[\frac{1}{\sqrt{2}}(1+j)]$	  & $-k $  & $ j$&  $i $ 	\\
			$[\frac{1}{2}(1-i-j+k)]$ & $ j$  & $ -k$&  $ -i$  & & 	$[\frac{1}{\sqrt{2}}(1-j)]$    & $k $  & $j $&  $-i $	\\

			$[\frac{1}{2}(1-i+j+k)]$ & $ -k$  & $ -i$& $ j$  & & 	$[\frac{1}{\sqrt{2}}(1+k)]$	  & $ j$  & $-i $&  $k $ \\	

				$[\frac{1}{2}(1+i-j+k)]$   & $ k$  & $ -i$&  $ -j$	& & 	$[\frac{1}{\sqrt{2}}(i-j)]$	   & $-j $  & $-i $&  $-k $	 \\
	
				 $[\frac{1}{\sqrt{2}}(1+i)]$	  & $ i$ & $ k$& $ -j$ 	& & 	&  & &				\end{tabular}
			\end{center} 
			\caption{Action of $\sfrac{\BO}{\mathbb{Z}_2}$ on $\BD_4$.}
 \label{tabla3.1}
			\end{table}

(4) The tetrahedral group $\T$ is isomorphic to the alternating group $\Alt(4)$, which has automorphism group $\Sym(4)$, acting by conjugation on the normal
subgroup  $\Alt(4)$. This corresponds to the action of the octahedral group $\OR\cong \Sym(4)$ on its normal subgroup $\T$, which is induced by the action of 
$\BO$ on the normal subgroup $\BT$.  In fact, it can be derived from Table \ref{tabla3.1} that the image of $\OR=\sfrac{\BO}{\mathbb{Z}_2}$ in $\Aut(\T)=\Aut(\sfrac{\BT}{\mathbb{Z}_2})$ is isomorphic to $\Sym(4)$. 

(5) Every automorphism of $\OR$ is inner.  

(6) The icosahedral group $\I$ is isomorphic to $\Alt(5)$, which is generated by $s=(12)(34)$ and $t=(135)$. Observe these generators satisfy $s^2=t^3=(st)^5=(1)$. The automorphism group of $\Alt(5)$ is $\Sym(5)$, whilst the outer automorphism group is isomorphic to $\mathbb{Z}_2$. In terms of permutations, the latter is generated by conjugation with an odd permutation, say $(35)$. This sends the generators listed above to $(12)(45)$ and $(153)$ respectively. The action of this automorphism $\varphi$ on conjugacy classes is described below. 

\begin{minipage}{0.45\textwidth}
\begin{flushleft} 
	\begin{table}[H] 
\begin{center}
  \begin{tabular}{ ||c | c  || }
    \hline
		Representative			 & Size 	 \\ \hline
			$(1)$     &  $1$\\
			$(123)$   &  $20$  \\
			$(12345)$ &  $12$		 \\ 
			$(13452)$		 & $12$  			 		 \\ 
			$(12)(34)$		 & $15$			  \\ \hline
			
\end{tabular}
\end{center}
\caption{Conjugacy classes. }
\label{conjAlt51}
\end{table}
\end{flushleft}
\end{minipage}
\begin{minipage}{0.48\textwidth}
\begin{flushright} 
\begin{table}[H]
\begin{center}
  \begin{tabular}{ ||c | c || }
    \hline		  
		$\mathcal{C}(xyzvw)$ 		 & $\varphi(\mathcal{C}(xyzvw))$  \\ \hline
		$\mathcal{C}(1)$		 & $\mathcal{C}(1)$  \\ 
		$\mathcal{C}(123)$	 	 & $\mathcal{C}(123)$		\\ 
		$\mathcal{C}(12345)$ 	 & $\mathcal{C}(13452)$	 		 			 \\ 
		$\mathcal{C}(13452)$ 	 & $\mathcal{C}(12345)$ \\ 
		$\mathcal{C}(12)(34)$ 	 & $\mathcal{C}(12)(34)$\\ \hline
\end{tabular}
\end{center}
\caption{Action of $\varphi$.}
\label{conjAlt5}
\end{table}
\end{flushright}
\end{minipage}   

\ \\
Where $\mathcal{C}(xyzvw)$ is the conjugacy class of a permutation $(xyzvw)\in \Alt(5)$. 

(7) The outer automorphism group of $\BT \subset \SU(2)$ is generated by an involution that exchanges the generators $s=\frac{1}{2}(1+i)(1+j), t=\frac{1}{2}(1+j)(1+i)$, which satisfy the relations $s^3=t^3=(st)^3$. This automorphism is given by conjugation with $\frac{1+j}{\sqrt{2}}\in  \BO \subset \SU(2)$.

(8) The outer automorphism group of $\BO \subset \SU(2)$ is generated by an involution $\varphi$ fixing $s$ and sending $t$ to $-t$, where $s= \frac{1}{2}(1+i+j+k) $ and $t= e^{\frac{i\pi}{4}}$ generate $\BO$.

(9) The outer automorphism group of $\BI \subset \SU(2)$ is generated by an involution $\psi$ which fixes $s$ and sends $t$ to $\frac{- \phi^{-1}-\phi i +k}{2}$ , where $s=\frac{1}{2}(1+i+j+k)$ and $t=\frac{\phi + \phi ^{-1}i+ j}{2}$  generate $\BI$. 

The action of the automorphisms $\varphi \in \Out(\BO)$ and $\psi \in \Out(\BI)$ on conjugacy classes $\mathcal{C}(x)$, for $x \in \BO$ or $x\in  \BI$ respectively, is described in the following tables.

\begin{minipage}{0.48\textwidth}
\begin{flushright} 
\begin{table}[H]
\begin{center}
  \begin{tabular}{ ||c | c | c|| }
    \hline
		Representative 			 & Size & Real parts	 \\ \hline
			$1$ & $1$ & $1$\\ 
			$-1$ & $1$ & $-1$			\\ 
			$s$ 		 &  	$8$ & $\frac{1}{2}$		 			 \\ 
			$t$ 		 &  $6$  & $\frac{1}{\sqrt{2}}$\\ 
			$s^2$ & $8$  & $-\frac{1}{2}$\\ 
			$t^2$ & $8$ & $0$			 \\ 
			$t^3$ 		 & $6$ &$- \frac{1}{\sqrt{2}}$ 			 			 \\ 
			$st$ 		 &  	$12$	 & $0$ \\ \hline
\end{tabular}
\end{center}
\caption{Conjugacy classes in  $\BO$ .}
\label{tabla7}
\end{table}
\end{flushright}
\end{minipage}
\begin{minipage}{0.45\textwidth}
\begin{flushleft} 
\begin{table}[H]
\begin{center}
  \begin{tabular}{ ||c | c | c|| }
    \hline		  
		$\mathcal{C}(x)$ 			 & $\varphi(\mathcal{C}(x))$ & $\Real(\varphi(x))$ 	 \\ \hline
			$\mathcal{C}(1)$ & $\mathcal{C}(1)$ & 1 \\ 
			$\mathcal{C}(-1)$ & $\mathcal{C}(-1)$ & 	-1		\\ 
			$\mathcal{C}(s)$ 		 &  	$\mathcal{C}(s)$ & $\frac{1}{2}$ 		 			 \\ 
			$\mathcal{C}(t)$ 		 & $\mathcal{C}(t^3)$   & $- \frac{1}{\sqrt{2}}$ \\ 
			$\mathcal{C}(s^2)$ & $\mathcal{C}(s^2)$  & $- \frac{1}{2}$ \\
			$\mathcal{C}(t^2)$ & $\mathcal{C}(t^2)$ & $0$ \\
			$\mathcal{C}(t^3)$ & $\mathcal{C}(t)$ & $\frac{1}{\sqrt{2}}$ 			 \\ 
			$\mathcal{C}(st)$ 		 & $\mathcal{C}(st)$ 		 & $0$ \\ \hline
\end{tabular}
\end{center}
\caption{Action of $\varphi$.}
\label{conj2O}
\end{table}
\end{flushleft}
\end{minipage}   

\begin{minipage}{0.48\textwidth}
\begin{flushright} 
	\begin{table}[H] 
\begin{center}
  \begin{tabular}{ ||c | c |c || }
    \hline
		Representative			 & Size 	&   Real parts \\ \hline
			$1$ & $1$ &   $1$ \\ 
			$-1$ & $1$ 			& 	$ -1$ \\ 
			$t$ 		 &  	$12$		 			&	  $\frac{1+\sqrt{5}}{4}$ \\ 
			$t^2$ 		 &  $12$			 			&   $-\frac{1-\sqrt{5}}{4}$  \\ 
			$t^3$ & $12$ &   $\frac{1-\sqrt{5}}{4}$ \\ 
			$t^4$ & $12$ 			& 	$-\frac{1+\sqrt{5}}{4} $ \\ 
			$s$ 		 & $20$ 			 			&	  $\frac{1}{2}$ \\ 
			$s^4$ 		 &  $20$			 			&   $-\frac{1}{2}$  \\ 
			$st$ 		 &  	$30$		 			&   $0$  \\ \hline
\end{tabular}
\end{center}
\caption{Conjugacy classes in  $\BI$. }
\label{tabla6}
\end{table}
\end{flushright}
\end{minipage}
\begin{minipage}{0.45\textwidth}
\begin{flushleft} 
\begin{table}[H]
\begin{center}
  \begin{tabular}{ ||c | c | c|| }
    \hline		   
		$\mathcal{C}(x)$ 			 & $\psi(\mathcal{C}(x))$ & $\Real(\psi(x))$ 	 \\ \hline
			$\mathcal{C}(1)$ & $\mathcal{C}(1)$ & 1 \\ 
			$\mathcal{C}(-1)$ & $\mathcal{C}(-1)$ & -1			\\ 
			$\mathcal{C}(t)$ 		 &  	$\mathcal{C}(t^3)$ & $\frac{1-\sqrt{5}}{4}$ 		 			 \\ 
			$\mathcal{C}(t^2)$ 		 & $\mathcal{C}(t^4)$   & $-\frac{1 + \sqrt{5}}{4}$ \\ 
			$\mathcal{C}(t^3)$ & $\mathcal{C}(t)$  & $\frac{1+ \sqrt{5} }{4}$ \\
			$\mathcal{C}(t^4)$ & $\mathcal{C}(t^2)$ & $- \frac{1- \sqrt{5}}{4}$ \\
			$\mathcal{C}(s)$ & $\mathcal{C}(s)$ & $\frac{1}{2}$ 			 \\ 
			$\mathcal{C}(s^4)$ 		 & $\mathcal{C}(s^4)$ & 	$- \frac{1}{2}$		 			 \\ 
			$\mathcal{C}(st)$ 		 & $\mathcal{C}(st)$ 		 & $0$ \\ \hline
\end{tabular}
\end{center}
\caption{Action of $\psi $. }
\label{conj2I}
\end{table}
\end{flushleft}
\end{minipage} 
\subsection{Freely acting splittable groups}
We regard splittable groups $A_1 \times A_2 \times A_3 \subset G $  with non-trivial factors $A_i\subset \SU(2) $, for $i=1,2,3$. The following technical result shall be used in the proof of Theorem \ref{main}.
\begin{Lem}\label{tech.prod}  Let $n,m \geq 2$ be integers, $k= \gcd (m,n)$ and $m_1=\frac{m}{k}, n_1= \frac{n}{k}$. Then, the solution set of 
$$\cos\left(\frac{2 \pi x}{n} \right)=\cos\left(\frac{2 \pi y}{m} \right) \ , \ (x,y)\in \mathbb{Z}^2, $$
is given by $\{ (nq+\varepsilon \ell n_1,\ell m_1) : q,\ell \in \mathbb{Z}, \varepsilon = \pm 1\}$.
\end{Lem} 
\pf Let $(x,y)\in \mathbb{Z}^2$ be a solution of either of the following equations
\begin{align}
mx = ny \mod  mn & & mx =  -ny \mod  mn.  \nonumber
\end{align}
say the first. In particular, there is an integer $ q \in \mathbb{Z}$, so that $m_1x-n_1y=mn_1 q$. Theorem \ref{lind} tells us that $(x,y)\in ((nq,0) + \mathbb{Z}(n_1,m_1))$. The converse of the assertion is straightforward.
\epf
\begin{Th} \label{main}
Any splittable subgroup acting freely on $M$, is up to permutation of the three factors and conjugation in $G= \SU(2)^3$, one in the following list.
\begin{table}[H]
\begin{center}
  \begin{tabular}{ ||c | c || }
    \hline
		Group 			 &   Conditions \\ \hline 
				$\mathbb{Z}_n \times \BI \times \BI$	 &   $ 2,3,5 \nmid  n$  \\ \hline 
					$\mathbb{Z}_n \times \BO \times \BI$	 & 	 \\
					$\mathbb{Z}_n \times \BO \times \BO$	 &	   \\ 
					$\mathbb{Z}_n \times \BT  \times \BO$	 &  $  2,3 \nmid  n$   \\ 
					$\mathbb{Z}_n \times \BT  \times \BI$	 &      \\
					$\mathbb{Z}_n \times \BT  \times \BT $	 &     \\ \hline 
					$\mathbb{Z}_n \times \mathbb{Z}_m \times  \mathbb{Z}_l$ &   $\gcd (n,m,l)=1$    \\ \hline 
					$\mathbb{Z}_n \times \mathbb{Z}_m \times \BD_{2l}$			 & $\gcd (n,m,2l)=1$  \\	 \hline
			    	$\mathbb{Z}_n \times \BD_{2m} \times  \BD_{2l}$&  $\gcd (n,2m,2l)=1$ 	 \\ \hline 
					$\mathbb{Z}_n \times \mathbb{Z}_m \times \BT $ &	   \\ 
					$\mathbb{Z}_n \times \mathbb{Z}_m \times \BO$ &  $2,3 \nmid \gcd (n,m)$	   \\ \hline 
					$\mathbb{Z}_n \times \BD_{2m} \times \BT $ & 	   \\
					$\mathbb{Z}_n \times \BD_{2m} \times \BO$ & $2, 3 \nmid \gcd (n,2m)$ \\ 
					$\mathbb{Z}_n \times \BD_{2m} \times \BI$ &    \\\hline	
					$\mathbb{Z}_n \times \mathbb{Z}_m \times \BI$	 & $2,3,5 \nmid  \gcd (n,m) $    \\ \hline
\end{tabular}
\end{center}
\caption{Freely acting splittable groups}
\label{splittable_table}
 
\end{table}
\end{Th}
\pf Observe first that a freely acting group  $A_1 \times A_2 \times A_3 \subset  G$ must have a cyclic factor  $\mathbb{Z}_n$ of odd order $n$, simply because otherwise $(-1-1-1)\in A_1 \times A_2 \times A_3 \subset  G$, see Table \ref{tabla4}, situation which is not allowed by Lemma \ref{crit}.
 A  splittable group $\mathbb{Z}_n \times A_2 \times A_3 \subset G$ acts freely on $M$ precisely if $\Real(A_2) \cap \Real(A_3) \cap \Real(\mathbb{Z}_n)=\{1\}$. However, since $n$ is odd, we always have 
 \[ \{-1,0, \pm \frac{1}{\sqrt{2}}, \frac{1}{2}, \frac{1}{4}(1\pm \sqrt{5})\} \cap \Real(\mathbb{Z}_n)= \emptyset ,\] 
 see Table \ref{tabla2} and \ref{tabla3}. Moreover, we see that  under the additional divisibility assumptions, which can be 
 read off from these tables, also the numbers $-\frac12$,  $\frac{1}{4}(-1\pm \sqrt{5})$ are not contained in $\Real(\mathbb{Z}_n)$. 
This observation gives us already all the freely acting splittable groups $\mathbb{Z}_n \times A_2 \times A_3 \subset G$ with $A_2, A_3 \subset \SU(2)$ of type E
satisfying the conditions displayed in Table \ref{splittable_table}. 

As for the cases involving at least one factor $A_2$ or $A_3 \subset \SU(2)$ of type AD, we observe that $ \Real (\mathbb{Z}_{n}) \cap \Real( \BD_{2m}) =\Real (\mathbb{Z}_{\gcd (n,2m)}), \Real (\mathbb{Z}_{n}) \cap \Real(\mathbb{Z}_m)= \Real (\mathbb{Z}_{\gcd (n,m)})$, which is a consequence of Lemma \ref{tech.prod}. Similar conclusions as before hold.
\epf
\section{Semi-splittable groups}\label{mainsection} 
We start the classification of freely acting semi-splittable subgroups $C \times D \subset G$ by considering the classes of such subgroups for which the quintuple $\mathcal{Q}(C)=\{A,A_0,B, B_0,\theta\}$ defining $C \subset  \SU(2)^2$ is so that $A,B \subset \SU(2)$ are non-trivial ADE groups and one of the following conditions holds. 
\begin{itemize}
\item[(a)] $A_0=A$, $B_0=B$.
\item[(b)] $A_0=B_0=\{1\}$ .
\item[(c)] $A_0 \neq\{1\}$ and $B_0 = \{1\}$.
\end{itemize}
Note that condition (a) is equivalent to $C \times D \subset G $ being splittable. Since splittable groups acting freely on $M=S^3\times S^3$ have already been classified
in the previous sections, we will assume from now on that our group is strictly semi-splittable. Observe also that under the condition (b), the group $C$ is conjugate
to the graph $\Gamma(\varphi,A)$ of some automorphism $\varphi$ of $A$. 
\begin{Def}
 A strictly semi-splittable group is said to be of type I or II if it satisfies condition (b) or (c) respectively, and of type  III otherwise.
 \end{Def}
\subsection{Type I groups}
We regard groups $\Gamma(\varphi,A) \times B \subset G$, where $\Gamma(\varphi,A)$ is the graph of some automorphism $\varphi: A \rightarrow A$ of an ADE group $A \subset  \SU(2)$ and $B\subset  \SU(2)$ is some non trivial ADE group. The following technical result shall be used recurrently in the forthcoming sections. 
\begin{Lem}\label{techn} Let $n \geq 2, r \in \mathbb{Z}_{n}^{\times}$, and set $k=\gcd (1+r,n), n_1=\frac{n}{k}, r_1= \frac{1+r}{k}$. Then, $n_1 \mathbb{Z}$ is the solution set of the equation: $x = -rx \mod  n $.
\end{Lem}
\pf Let $x \in \mathbb{Z}$ be a solution of the equation in question. There is in particular an integer $y \in \mathbb{Z}$, so that $(1+r)x+ny=0$. If we consider the latter as an equation in $(x,y)\in \mathbb{Z}^2$, then $(x,y)=t(n_1,-r_1)$ for some $t \in \mathbb{Z}$, see Theorem \ref{lind}. This shows the claim as any number $tn_1\in \mathbb{Z}$ is a solution of: $x = -rx \mod  n $.  
\epf
\begin{Th}\label{typeB}
Let $\Gamma(\varphi,A) \times B \subset G$ be a freely acting type I group which is not a subgroup of a freely acting 
splittable group. Then, it is either $\Gamma(\varphi,\BI)\times \mathbb{Z}_n$ with $3 \nmid n$, or one of the following groups. 

\begin{table}[H]
\begin{center}
\begin{tabular}{ ||c|c|c| c|| }
\hline
\multicolumn{4}{ ||c ||}{$\Gamma(\tau_{a,b},\BD_{2n})  \times B \ , \ n>2, \ a \neq \pm 1 \mod 2n$ }  \\
\hline
$\tilde{k}_1>1$& $\tilde{k}_2>1$  & $B$ & Conditions \\ \hline
$yes$& $yes$ & $\mathbb{Z}_m$ & $4 \nmid m \ \wedge \ \gcd (\tilde{k}_1,m)=\gcd (\tilde{k}_2,m)=1 $ \\ \hline
$yes$& $no$		& $\mathbb{Z}_m$ & $4 \nmid m \ \wedge \ \gcd (\tilde{k}_1,m)=1\ $ \\ \hline
$no$& $yes$		& $\mathbb{Z}_m$ & $4 \nmid m \ \wedge \ \gcd (\tilde{k}_2,m)=1\ $ \\ \hline
$no$&$no$ 			& $\mathbb{Z}_m$&  $4\nmid m$ \\ \hline 
\end{tabular}
\end{center}
\end{table}
\begin{table}[H]
\begin{center}
\begin{tabular}{ ||c|c|c| c|| }
\hline
\multicolumn{4}{ ||c| }{$\Gamma(\varphi(r),\mathbb{Z}_n) \times B, \ \ r \neq \pm 1 \mod n$ }  \\
\hline
 $k_1>1$& $k_2>1$  & $B$ & Conditions \\ \hline
 & & $\mathbb{Z}_m$& $\gcd (m,k_1)= \gcd (m,k_2)=1$\\
 $yes$& $yes$ & $\BD_{2m}$ & $ \gcd (2m,k_2)=\gcd (2m,k_1)=1$ \\
&   & $ \BT,\BO$  & $2, 3 \nmid  k_1,k_2$ \\
 & & $\BI$& $2, 3,5 \nmid  k_1,k_2$ \\ \hline
 
 &  &$\mathbb{Z}_m$ & $\gcd (k_1,m)=1$ \\ 
 $yes$& $no$		& $\BD_{2m}$ & $ \gcd (k_1,2m)=1$ \\ 
 & 	 & $\BT,\BO$& $2,3 \nmid  k_1$ \\ 
& & $\BI$&  $2,3,5 \nmid  k_1$ \\ \hline
&  &$\mathbb{Z}_m$ & $\gcd (k_2,m)=1$ \\ 
 $no$& $yes$		& $\BD_{2m}$ & $\gcd (k_2,2m)=1$ \\ 
 & 	 & $\BT,\BO$& $2,3 \nmid  k_2$ \\ 
& & $\BI$&  $2,3,5 \nmid  k_1$ \\ \hline
$no$&$no$ 			& All&  - \\ \hline 
\end{tabular}
\end{center}
\caption{Type I freely acting subgroups. }
\end{table}
Where $k_1=\gcd (1+r,n)$, $k_2=\gcd (1-r,n)$, $\tilde{k}_1=\gcd (1+a,2n) $ and $\tilde{k}_2=\gcd (1-a,2n)$.
\end{Th}
\pf Let $A \subset \SU(2)$ be an ADE group and define $$ \mathcal{W}(\varphi, A)= \{ \Real(x)  :  \Real(x) = \Real ( \varphi(x) ) \}. $$ If $\varphi \in \Inn(A)$, then $\mathcal{W}(\varphi , A)=  \Real{(A)}$. In consequence, a group $\Gamma(\varphi,A) \times B \subset G$ with $\varphi \in \Inn(A)$ that acts freely on $M$ must be a subgroup of a freely acting splittable group. We consider therefore just outer automorphisms of $A \subset \SU(2)$. Observe that once we have calculated $\mathcal{W}(\varphi, A ) \subset \mathbb{R}$, the precise conditions to impose on $\Gamma(\varphi,A) \times B \subset G$ to make it act freely on $M$ are easily obtained from 
Tables \ref{tabla4}, \ref{tabla2} and \ref{tabla3} as in the case of splittable groups. For this reason, we just calculate $\mathcal{W}(\varphi, A )$ for any ADE group $A \subset \SU(2)$ and $\varphi \in \Out(A)$:

(a) let $A\subset  \SU(2)$ be a type E group. Since $\Out (A)$ is in this case isomorphic to $\mathbb{Z}_2$, we must consider a single automorphism of the group $A \subset  \SU (2)$. Furthermore, as the only non-trivial class in $\Out(\BT)$ is represented by an automorphism which can be obtained by conjugation with some element in $\SU (2)$, we are left with $A \in \{\BO,\BI\}$. In these cases, we read from Tables \ref{tabla7}- \ref{conj2I} that 
$$ \mathcal{W}(\varphi , \BO )  = \mathcal{W}(\psi , \BI ) =\left\{ 1 ,-1,  \frac{1}{2}, - \frac{1}{2} ,  0\right\}.$$

(b) Let $A=\mathbb{Z}_n$ and $\varphi(r)$ be a non-trivial outer automorphism of $\mathbb{Z}_n$. Elements in $\mathcal{W}( \varphi(r), \mathbb{Z}_n)$ are readily seen to be obtained by solving 
\begin{align}
x = \pm rx \mod  n. \label{equ.e}
\end{align}
Further, we can suppose $r \neq \pm 1 \mod n$, as otherwise $\Real(\mathbb{Z}_n)=\mathcal{W}(\varphi(r), \mathbb{Z}_n)$. If $\gcd (1 \pm r,n)=1$, equation (\ref{equ.e}) admits only the trivial solution and, thus, any group $\Gamma(\varphi(r),\mathbb{Z}_n) \times B \subset  G$ acts freely on $M$. Now, suppose that $k_1=\gcd (1+r,n)>1,\ k_2=\gcd (1-r,n)=1$ and write $n=k_1n_1, \ 1+r=k_1r_1$. The solution set of equation (\ref{equ.e}) is $n_1 \mathbb{Z}$, see Lemma \ref{techn}, and so
	\begin{align}
	\mathcal{W}( \varphi(r) , \mathbb{Z}_n )&= \left\{ \cos\left(\frac{2\pi z}{k_1} \right) \ : \ z \in \mathbb{Z}\right\}. \nonumber
	\end{align}
	 If $k_1 ,  k_2 >1$, we find out that
	\begin{align}
	\mathcal{W}(\mathbb{Z}_n, \varphi(r))&= \left\{ \cos\left(\frac{2\pi z}{k_i} \right) \ : \ z \in \mathbb{Z},\ i=1,2 \right\}. \nonumber
	\end{align}

(c) At last, consider $A=\BD_{2n}$. Since automorphism of $\BD_4$ are obtained by conjugation with an element in $\BO \subset \SU(2)$, we can suppose that $n>2$. Let $\tau_{a,b}$ be a non-trivial outer automorphism of $\BD_{2n}$. Since 
	\begin{align}
\Gamma(\tau_{a,b},\BD_{2n})&= \left\{ \left(e^{\frac{i \pi x}{n}}, e^{\frac{i \pi a x}{n}} \right), \left( je^{\frac{i \pi x}{n}},j e^{\frac{i \pi (a (x-1)+b+1)}{n}} \right) \ : \ x \in \mathbb{Z} \right\}, \nonumber
\end{align}
an element in $\mathcal{W}(\BD_{2n},\tau_{a,b})$ is either zero or it is obtained by solving the equation: $x = \pm ax \mod  2n$. \epf
\subsection{Type II groups}
We regard groups $C \times D \subset  G$ so that the quintuple $\mathcal{Q}(C)=\{A,A_0,B,B_0,\theta\}$ defining $C\subset \SU(2)^2$ is so that $B_0=\{ 1\}$ and 
$A \subset \SU(2)$ admits a quotient isomorphic (via $\theta$) to some ADE group $B \subset \SU(2)$. As explained in the remark on page \pageref{Remark}, we shall identify here $A/A_0$
with $F=B$ and consider 
$\theta$ as an automorphism of $F$.
Tables \ref{tabla-p} and \ref{tabla-p.1} (which are borrowed from \cite{FIG}) display all non-trivial normal subgroups $A_0\subset \SU(2)$ of an ADE group $A \subset \SU(2)$ and the isomorphism type of the corresponding quotient. 

\begin{minipage}{0.5\textwidth}
\begin{flushleft} 
\begin{table}[H]
\begin{center}
  \begin{tabular}{ ||c | c || }
    \hline
		$A_0 \triangleleft A$ 			 & $A/A_0$	 \\ \hline
			$\mathbb{Z}_k \triangleleft \mathbb{Z}_{kl} $ & $ \mathbb{Z}_l $     \\ 
			$\mathbb{Z}_{2k} \triangleleft \BD_{2kl}$ & $\D_{2l}$ 	 \\  
			$\mathbb{Z}_{2k+1} \triangleleft \BD_{2l(2k+1)}$		 & $\BD_{2l}$  			 	 \\ 
			$\mathbb{Z}_{2k+1} \triangleleft \BD_{2(2k+1)}$		 & $\mathbb{Z}_4$  \\
			$\BD_{2k} \triangleleft \BD_{4k}$		 & $\mathbb{Z}_2$  			  \\ 
			$\mathbb{Z}_2 \triangleleft \BT$  & $\T$  \\ \hline
			
\end{tabular}
\end{center}
\caption{Subgroups  I. }
\label{tabla-p}
\end{table}
\end{flushleft}
\end{minipage}
\begin{minipage}{0.5\textwidth}
\begin{flushright} 
\begin{table}[H]
\begin{center}
  \begin{tabular}{ ||c | c || }
    \hline
		$A_0 \triangleleft A$ 			 & $A/A_0$ 	 \\ \hline
		
$\BD_4 \triangleleft \BT$	 & $\mathbb{Z}_3$      \\ 
		$\mathbb{Z}_2 \triangleleft \BO$ 	 &  $\OR$ 	 \\ 
			$\BD_4 \triangleleft \BO$ 		 &  $\D_6$			 	 \\ 
			$\BT\triangleleft \BO $ 		 &  $\mathbb{Z}_2$			  \\ 
			$\mathbb{Z}_2 \triangleleft \BI$  & $\I$\\ \hline
\end{tabular}
\end{center}
\caption{Subgroups II.}
\label{tabla-p.1}
\end{table}
\end{flushright}
\end{minipage}

\begin{cav}
It should be stressed here that according to our convention, the group $\BD_{2n}$ is defined for $n > 1$.
\end{cav} 
\begin{Th}\label{type3} 
Any type II freely acting group $C \times D \subset G$ that is not a subgroup of a freely acting splittable group belongs to the following list 
\begin{table}[H]
\begin{center}
\begin{tabular}{ ||c|c|| }
\hline
 Group & Conditions \\ \hline
 $\mathcal{G}(\BT,\BD_4,\mathbb{Z}_3 , \{1\}, \varphi(r)) \times \mathbb{Z}_n$ & $3 \nmid  n$ \\ \hline
 $\mathcal{G}(\BT,\BD_4,\mathbb{Z}_3 , \{1\}, \varphi(r)) \times \BD_{2l}$ & $3 \nmid  2l$ \\ \hline
 $\mathcal{G}(\BD_{4k},\BD_{2k}, \mathbb{Z}_2, \{1\}, \ID) \times D \ $ & \\ 
$ \mathcal{G}(\BO,\BT, \mathbb{Z}_2, \{1\}, \ID) \times D$ & $-$ \\ 
$ \mathcal{G}(\BD_{2k},\mathbb{Z}_{2k},\mathbb{Z}_2, \{1\}, \ID) \times D$ & \\  \hline 
\end{tabular}
\end{center}
\end{table}

 \begin{table}[H]
\begin{center}
\begin{tabular}{ ||c|c|c| c|| }
\hline
\multicolumn{4}{ ||c| }{$\mathcal{G}( \mathbb{Z}_{kl}, \mathbb{Z}_l, \mathbb{Z}_k, \{1\}, \varphi(r)) \times D, \ lr \neq \pm 1 \mod lk$ }  \\
\hline
 $k_1>1$& $k_2>1$  & $D$ & Conditions \\ \hline
 & & $\mathbb{Z}_m$& $\gcd (m,k_1)=\gcd (m,k_2)=1$\\
 $yes$& $yes$ & $\BD_{2m}$ & $ \gcd (2m,k_2)=\gcd (2m,k_1)=1$ \\
&   & $ \BT,\BO$  & $2, 3 \nmid  k_1,k_2$ \\
 & & $\BI$& $2, 3,5 \nmid  k_1,k_2$ \\ \hline 
 &  &$\mathbb{Z}_m$ & $\gcd (k_1,m)=1$ \\ 
 $yes$& $no$		& $\BD_{2m}$ & $ \gcd (k_1,2m)=1$ \\ 
 & 	 & $\BT,\BO$& $2,3 \nmid  k_1$ \\ 
& & $\BI$&  $2,3,5 \nmid  k_1$ \\ \hline
&  &$\mathbb{Z}_m$ & $\gcd (k_2,m)=1$ \\ 
 $no$& $yes$		& $\BD_{2m}$ & $\gcd  (k_2,2m)=1$ \\ 
 & 	 & $\BT,\BO$& $2,3 \nmid  k_2$ \\ 
& & $\BI$&  $2,3,5 \nmid  k_1$ \\ \hline
$no$&$no$ 			& All&  - \\ \hline 
\end{tabular}
\end{center}
\end{table}
\begin{table}[H]
\begin{center}
\begin{tabular}{ ||c|c|c| c|| }
\hline
\multicolumn{4}{ ||c ||}{$ \mathcal{G}(\BD_{2l(2k+1)},\mathbb{Z}_{2k+1},\BD_{2l}, \{1\}, c_{g} \circ \tau_{a,b})  \times D,  $ }  \\
\multicolumn{4}{ ||c ||}{$\  \ a(2k+1) \neq \pm1 \mod 2l(2k+1), \ l>2 $} \\
\hline
$\tilde{k}_1>1$& $\tilde{k}_2>1$  & $D$ & Conditions \\ \hline
$yes$& $yes$ & $\mathbb{Z}_m$ & $ \gcd (\tilde{k}_1,m)=\gcd (\tilde{k}_2,m)=1 \ \wedge 4\nmid m $ \\ \hline
$yes$& $no$		& $\mathbb{Z}_m$ & $\gcd (\tilde{k}_1,m)=1 \ \wedge 4 \nmid m $ \\ \hline
$no$& $yes$		& $\mathbb{Z}_m$ & $\gcd (\tilde{k}_2,m)=1 \ \wedge 4 \nmid m $ \\ \hline
$no$&$no$ 			& $\mathbb{Z}_m$& $4 \nmid m$  \\ \hline 
\end{tabular}
\end{center}
\caption{Type II freely acting subgroups. }
\end{table}
Where $k_1=\gcd (1+lr,kl)$, $k_2=\gcd (1-lr,kl)$, $\tilde{k}_1=\gcd (1-a(2k+1),2l(2k+1))$ and $\tilde{k}_2=\gcd (1+a(2k+1),2l(2k+1))$.
\end{Th}
\pf Let $C \times D \subset  G$ be a freely acting type III group, $\mathcal{Q}(C)=\{A,A_0,B,\{ 1 \},\theta \}$ be the quintuple defining $C \subset \SU(2)^2$ and define $\mathcal{W}(C) = \{ \Real (x)\ :  \Real (x)= \Real (y) \ , \ (x,y)\in C \}$. The following are the group triples $(A,A_0,B)$ that give rise to a valid choice for $\mathcal{Q}(C)$, 
that is such that the quotient $B=A/A_0$ is an ADE group, see Tables \ref{tabla-p} and \ref{tabla-p.1}.   
\begin{table}[H]
\begin{center}
  \begin{tabular}{ ||c | c  || }
    \hline
		\textit{B} 			 &   $(A,A_0)$ \\[0.1cm] \hline 
		$\mathbb{Z}_k$	&	$(\mathbb{Z}_{kl},\mathbb{Z}_l)$  \\ 
		$\mathbb{Z}_4$ & $(\BD_{2(2k+1)},\mathbb{Z}_{2k+1} )  $\\
		$\mathbb{Z}_3$	&		$(\BT,\BD_4)$		\\
		
		$\mathbb{Z}_2$	&	$(\BD_{4k},\BD_{2k})$   \\ 
		$\mathbb{Z}_2$	&	$(\BO,\BT)$   \\
		$ \mathbb{Z}_2$ & $(\BD_{2k},\mathbb{Z}_{2k})$ \\
		$\BD_{2l}$	&	$(\BD_{2l(2k+1)},\mathbb{Z}_{2k+1})$ \\[0.1cm] \hline
\end{tabular}
\end{center}
\caption{Triples $(A,A_0,B)$.}
\label{tabla9}
\end{table}

We calculate $\mathcal{W}(C)$ case by case for triples $(A,A_0,B)$ that are part of $\mathcal{Q}(C)$ in order to read off
the conditions for $\mathcal{W}(C) \cap \Real (D) = \{ 1\}$. For later use we recall, see equation (\ref{fiber-prod.1}), that 
the group $C$ is a fibered product associated with the maps $\alpha : A \rightarrow B$ and $\beta =  \ID_B: B\rightarrow B$. 

(a) Let $(A,A_0,B)= (\mathbb{Z}_{kl},\mathbb{Z}_l, \mathbb{Z}_k)$ and $\theta = \varphi (r)$ some automorphism of $\mathbb{Z}_k$. We have that 
\begin{align}
C&=\{ ([y]_{kl},[ry]_k) \ : \ y \in \mathbb{Z} \} , \nonumber
\end{align}
for some $r \in \mathbb{Z}_{k}^{\times}$. Lemma \ref{techn} helps us determining the solution set of 
\begin{align}
(1 \pm lr) y = 0 \ \mod   kl, \nonumber 
\end{align}
whenever $lr \neq \pm 1 \mod lk$. The latter condition can be assumed as otherwise $\mathcal{W}(C)=\Real(\mathbb{Z}_{kl})=\Real(\mathbb{Z}_{k})$ and 
$C\times D$ is then a subgroup of a freely acting splittable group $A\times B \times D$.  We see that under this condition 
$$\mathcal{W}(C)=\left\{ \cos\left( \frac{2 \pi x}{k_i } \right) \ : \ x \in \mathbb{Z},\  i=1,2 \right\}.$$

(b) Let $(A,A_0,B)= (\BD_{2(2k+1)}, \mathbb{Z}_{2k+1}, \mathbb{Z}_4)$. The map $\alpha : A \rightarrow B$ defining $\mathcal{G}( \BD_{2(2k+1)},  \mathbb{Z}_{2k+1}, \mathbb{Z}_4, \{1\}, \varphi(r)) \subset \SU(2)^2$ is given by 
\begin{align}
	\alpha &: z \mapsto \alpha (z)= \left\{ 
  \begin{array}{l l}
    \text{$1 $} & \quad \mbox{if} \quad \text{$z = e^{\frac{i \pi  x }{2k+1} },\quad  x = 0  \mod   2 $,  } \\
		\text{$e^{\frac{i \pi  r }{2} }  $} & \quad \mbox{if} \quad  \text{$z = je^{\frac{i \pi  x }{2k+1} }, \quad   x = 0  \mod   2$,} \\
		\text{$e^{i \pi r }  $}& \quad \mbox{if} \quad\text{$z = e^{\frac{i \pi  x }{2k+1} },\quad   x = 1 \mod  \ 2 $,  } \\
		\text{$ e^{\frac{3i \pi  r }{2} }$}& \quad \mbox{if} \quad \text{$z =j e^{\frac{i \pi  x }{2k+1} },\quad    x = 1 \mod   2$}.
  \end{array} \right.   \nonumber
	\end{align}
	It is not difficult to verify that $( e^{\frac{i \pi  x }{2k+1} }, e^{\pi i } ) \in C \subset  \SU(2)^2$, for $x$ odd and any choice of $r \in \mathbb{Z}_4^{\times}$, hence $-1 \in \mathcal{W}(C)$. Since $C\times D$ acts freely it follows that  $-1 \not\in \Real (D)$ and, hence, $D=\mathbb{Z}_n$ for $n$ odd. Using $\mathcal{W}(C)\subset 
	\Real (\mathbb{Z}_4)= \{ 0, \pm 1\}$, this easily implies that $C \times D \subset G$ is
	 a subgroup of a freely acting  splittable group $\BD_{2(2k+1)} \times \mathbb{Z}_4 \times \mathbb{Z}_n$.  
	 
(c) Let $(A,A_0,B)= (\BT,\BD_4, \mathbb{Z}_3)$. The map defining $\mathcal{G}(\BT,\BD_4, \mathbb{Z}_3, \{1\}, \\ \varphi(r)) \subset \SU(2)^2$ is now given as follows
\begin{align}
	\alpha &: z \mapsto \alpha (z)= \left\{ 
  \begin{array}{l l}
    \text{$1$} & \quad \mbox{if} \quad\text{$z \in \BD_4$,}\\
    \text{$e^{\frac{2\pi i r}{3} }$} & \quad \mbox{if} \quad \text{$z \in  \frac{-1+i+j+k}{2} \BD_4 $,  } \\
		\text{$e^{\frac{4 \pi i r }{3} }$} & \quad \mbox{if} \quad \text{$z \in  \frac{1+i+j+k}{2} \BD_4 $.}
  \end{array} \right.  \nonumber
	\end{align}
We can check that $\pm \frac{1}{2} \in \mathcal{W}(C)$ for any choice of $r \in \mathbb{Z}_3 ^{\times}$. It follows that $\mathcal{W}(C)=\left\{1, \pm \frac{1}{2} \right\} =\Real (\mathbb{Z}_3)$, and so no new freely acting group is obtained in this way.

(d) Let $(A,A_0,B)= (\BD_{4k},\BD_{2k}, \mathbb{Z}_2)$. We have 
\begin{align}
C&= \{ (e^{\frac{i \pi y }{k}},1),(je^{\frac{i \pi y}{k}},1),(e^{\frac{i \pi (2y+1)}{2k}},-1),(j e^{\frac{i \pi (2y+1)}{2k}},-1) \ : \ y \in \mathbb{Z} \} , \nonumber
\end{align} 
and hence $\mathcal{W}(C)=\{ 1\}$. We conclude that any group $\mathcal{G}(\BD_{4k},\BD_{2k},\mathbb{Z}_2  , \{1\}, \ID )  \times D \subset G$ acts freely on $M$. 

(e) Let $(A,A_0,B)= (\BO,\BT,\mathbb{Z}_2)$. In this case, 
\begin{align}
	\alpha &: z \mapsto \alpha (z)= \left\{ 
  \begin{array}{l l}
    \text{$1$} & \quad \mbox{if} \quad \text{$z \in \BT$, }\\
    \text{$-1$} & \quad \mbox{if} \quad \text{$z \in e^{\frac{i \pi}{4} }\BT $. } 	
  \end{array} \right. \nonumber
	\end{align}
	 Because $\Real(e^{\frac{i \pi }{4}}\BT)= \{ 0, \pm \frac{1}{\sqrt{2}}, \pm \frac{1}{2 \sqrt{2}} \}$, we have that $\mathcal{W}(C)=\{1\}$. Therefore, every subgroup $\mathcal{G}(\BO,\BT, \mathbb{Z}_2, \{1\}, \ID) \times D \subset G$ acts freely on $M$.
	 
(f) Let $(A,A_0,B)= (\BD_{2k},\mathbb{Z}_{2k},\mathbb{Z}_2)$. We have that
\begin{align}
	\alpha &: z \mapsto \alpha (z)= \left\{ 
  \begin{array}{l l}
    \text{$1$}& \quad \mbox{if} \quad \text{$z=e^{\frac{i \pi x}{k}}$,}\\
    \text{$-1$} & \quad \mbox{if} \quad\text{$z=je^{\frac{i \pi x}{k}} $.} 	
  \end{array} \right. \nonumber
	\end{align}
	From which, we see that $$C=\{ (e^{\frac{\pi i x}{k} },1) , (je^{\frac{\pi i x}{k} },-1) \ : \ x \in  \mathbb{Z} \},$$
and so $\mathcal{W}(C)=\{1\}$. Any subgroup $\mathcal{G}(\BD_{2k},\mathbb{Z}_{2k},\mathbb{Z}_2, \{1\}, \ID) \times D \subset G$ will then act freely on $M$.   

(g) Let $(A,A_0,B)= (\BD_{2l(2k+1)},\mathbb{Z}_{2k+1},\BD_{2l})$ and $l>2$. Since $\BD_{2l}$ is non-commutative, we must consider also inner automorphisms in the construction of our group $C \subset \SU(2)^2$. Denote by $c_{g(w)}$ the conjugation map by an element $g(w) \in \{e^{\frac{i \pi w}{l}} , j e^{\frac{i \pi w}{l} }  :  w \in \mathbb{Z}\} \subset  \BD_{2l}$. The map defining $C \subset \SU(2)^2$ is given by 
\begin{align}
	\alpha &: z \mapsto \alpha (z)= \left\{ 
  \begin{array}{l l}
    \text{$ c_{g(w)}(e^{\frac{i \pi x a}{l} }) $} & \quad \mbox{if} \quad \text{$ z=e^{ \frac{i \pi x }{ l(2k+1) } } $, }\\
    \text{$ c_{g(w)}(je^{\frac{i \pi (a(x-1)+b+1)}{l} })$} & \quad \mbox{if} \quad \text{$ z=je^{\frac{i \pi x }{ l(2k+1)}} $. } 	
  \end{array} \right. \nonumber
	\end{align}
 where $(a,b) \in \mathbb{Z}_{2l}^{\times} \times \mathbb{Z}_{2l}$. In particular, $C \subset \SU(2)^2 $ equals
\begin{align}
\{ (e^{\frac{i \pi x }{l(2k+1)}}, c_{g(w)}(e^{\frac{i \pi x a}{l} })),  (je^{\frac{i \pi y }{l(2k+1)}},  c_{g(w)}(je^{\frac{i \pi (a(y-1)+b+1)}{l} }) ) \ : \  x ,y\in \mathbb{Z} \}.   \nonumber 
\end{align} 
  An element in $\mathcal{W}(C)$ is thus either zero, or it is obtained from a solution of  
\begin{align}
(1 \pm a(2k+1))x = 0 \mod 2l(2k+1). \nonumber
\end{align}  
This situation that can be treated analogously as in case (a). At last, let $l=2$ and $c_{\tilde{w}}$ be the conjugation map by an element $\tilde{w} \in \BO$. We find that
\begin{align}
C=\{ (e^{\frac{i \pi x }{2(2k+1)}}, c_{\tilde{w}}(e^{\frac{i \pi x}{2}   })),  (je^{\frac{i \pi y }{2(2k+1)}}, c_{\tilde{w}}(je^{\frac{i \pi y}{2}   })) \ : \  x,y \in \mathbb{Z} \}  , \nonumber 
\end{align}
and so, $\mathcal{W}(C)=\Real( \mathbb{Z}_4)$. No new freely acting subgroups are thus obtained.
\epf
\subsection{Type III groups} 
We seek now to distinguish freely acting type  III groups  $C \times D \subset G$. To begin with, we assert that such groups must fulfill rather restrictive conditions a priori.  
\bp \label{form.strict}
Let $C \times D \subset G$ be a type  III freely acting group and $C=\mathcal{G}(A,A_0,B,B_0, \theta) \subset \SU(2)^2$. Then (up to interchanging 
the roles of $A$ and $B$)
\begin{itemize}
\item[(a)]  $D= \mathbb{Z}_{2k+1}$ or 
\item[(b)] $(A,A_0,B,B_0)\in \{( \BD_{2l(2k+1)}, \mathbb{Z}_{2k+1}, \BD_{2l(2p+1)} , \mathbb{Z}_{2p+1}),  ( \mathbb{Z}_{3(2k+1)}, \mathbb{Z}_{2k+1},\BT$, $\BD_4)$, 
$(\mathbb{Z}_{(2k+1)l}, \mathbb{Z}_{2k+1}, \mathbb{Z}_{pl}  , \mathbb{Z}_p)$ , 
$(\mathbb{Z}_{2(2k+1)}, \mathbb{Z}_{2k+1}, \BO, \BT )$, 
$(\mathbb{Z}_{2(2k+1)}, \mathbb{Z}_{2k+1}$, $\BD_{4p}, \BD_{2p})$, 
$(\mathbb{Z}_{4(2k+1)}, \mathbb{Z}_{2k+1}, \BD_{2(2p+1)}, \mathbb{Z}_{2p+1})$, 
$( \BD_{2(2k+1)}, \mathbb{Z}_{2k+1},  \BD_{2(2p+1)}$, $\mathbb{Z}_{2p+1})$, $(\BD_{2(2k+1)}, \mathbb{Z}_{2k+1}, \mathbb{Z}_{4p}, \mathbb{Z}_p) \}.$
\end{itemize}
\ep
\pf  Let $C \times D \subset G$ be a type  III freely acting subgroup and $\mathcal{Q}(C)=\{A,A_0,B,B_0,\theta \}$ be the quintuple defining $C \subset  \SU(2)^2$ via the homomorphisms $\alpha$ and $\beta$   in equation (\ref{fiber-prod.1}). If $D \subset \SU(2)$ is a DE group or a cyclic group of even order, and the groups $A,B \subset \SU(2)$ belong to the following list:  groups of type E, $ \mathbb{Z}_n$ and $\BD_{2m}$, where $n$ and $m$ are powers of $2$.  A short glance at Tables \ref{tabla4}, \ref{tabla-p} and \ref{tabla-p.1} reveals that $-1 \in A_0 \cap B_0 \cap D$. This can not happen if $C \times D \subset G$ acts freely on $M$, as otherwise $(-1,-1,-1)\in C \times D \subset G$. By listing the remaining possibilities according to Tables \ref{tabla-p}, \ref{tabla-p.1} and Theorem \ref{GOURSAT}, we encounter the necessity of fulfilling (at least) one of the conditions in the theorem  if the group $C\times D$ acts freely.
\epf
 The technical result below will help us distinguishing type  III freely acting groups of the forms stated in Proposition \ref{form.strict}.
\begin{Lem}\label{res}  Let $m,n,p \geq 2$. The equation, 
$$\cos\left( \frac{2 \pi x}{3n}\right) = c, $$
has a solution
\begin{itemize}
\item[(a)] $x \in 1+3 \mathbb{Z}$ for $c=\frac{1}{2}$, precisely if $6 | (n-2)$ or $6|(n+2)$.
\item[(b)] $x \in 1+3 \mathbb{Z}$ for $c=-\frac{1}{2}$, precisely if $3 | (n-1)$ or $3 | (n+1)$.
\item[(c)] $x \in 2+3 \mathbb{Z}$ for $c=\frac{1}{2}$, precisely if $6 | (n-4)$ or $6 |(n+4)$.
\item[(d)] $x \in 2+3 \mathbb{Z}$ for $c=-\frac{1}{2}$, precisely if $3 | (n-2)$ or $3 | (n+2)$.
\end{itemize}
\end{Lem}
\bt \label{q.final}  For the groups $C= \mathcal{G}(A,A_0,B,B_0,\theta)$ considered below, which include the ones in part (b) of Proposition \ref{form.strict}, we describe in each case all type III groups $C \times D \subset G$ that acts freely on $M$.

(a)  Let $C= \mathcal{G}(\mathbb{Z}_{kl}, \mathbb{Z}_{k}, \mathbb{Z}_{pl}  , \mathbb{Z}_p,  \varphi(r))$, then $C \times D$ belongs to the following list. 
\begin{table}[H]
\begin{center}
\begin{tabular}{ ||c|c|| }
\hline
\multicolumn{2}{ ||c|| }{$  \mathcal{G}(\mathbb{Z}_{kl}, \mathbb{Z}_{k}, \mathbb{Z}_{pl}  , \mathbb{Z}_p,  \varphi(r))  \times D$ }   \\ \hline
 $D$ &  Conditions  \\ \hline 
 $\mathbb{Z}_{n}$& $\gcd (n, ms)=1$ \\ \hline  
 $\BD_{2n}$&  $ \gcd (2n,ms)=1$ \\ \hline  	  
 $\BT,\BO$& $2,3 \nmid  ms$ \\ \hline
  $\BI$	 & $2,3,5 \nmid  ms$ \\ \hline
\end{tabular}
\end{center}
\end{table}
Where the above conditions are required for both values of $$m=m_\varepsilon=\gcd (p -\varepsilon  kr, kl,lpk)),$$ for $\varepsilon = \pm 1$, and $$ s =\frac{kl}{\gcd (kl,lpk)}.$$

(b) Let $C = \mathcal{G}(\mathbb{Z}_{3n}, \mathbb{Z}_n,\BT,\BD_4, \varphi(r)) $ and assume that $C \times D$ is not a subgroup of a splittable freely acting group. Then it is either $\mathcal{G}(\mathbb{Z}_{3n}, \mathbb{Z}_n, \BT,\BD_4,$ $ \varphi(r))  \times \mathbb{Z}_{2p+1}$  with $n$ even, subject to the conditions:  
$$3 \nmid (n-1),(n+1),(n+2),(n-2) ,$$
or one of the following 
\begin{table}[H]
\begin{center}
\begin{tabular}{ ||c|c|c|c|| }
\hline
\multicolumn{4}{ ||c|| }{$ \mathcal{G}(\mathbb{Z}_{3n}, \mathbb{Z}_n,\BT,\BD_4, \varphi(r)) \times D$, \quad  $2 \nmid n$ }   \\ \hline
 $ 6 | (n \pm 2) \vee \ 6 | (n\pm 4) $  & $ 3 | (n \pm 1) \ \vee \ 3 | (n\pm 2) $ & $D$ & Conditions  \\ \hline
$yes$ & $yes$ & $\mathbb{Z}_m$ 	  & $3 \nmid m$ 	\\ 
& & $\BD_{2m}$ & $3 \nmid  m$ \\ \hline
 $yes$ & $no$ & $\mathbb{Z}_m$ 	  & $6 \nmid m$ 	\\ 
 & & $\BD_{2m}$ & $3 \nmid  m$ \\ \hline
 $no$ & $yes$ & $\mathbb{Z}_m$ 	  & $3 \nmid m$ 	\\ 
& & $\BD_{2m}$ & $3 \nmid  m$ \\ \hline
 $no$  &	$no$ 	 &	All 	 & -\\ \hline
\end{tabular}
\end{center}
\caption{Type III freely acting subgroups.}
\label{stt.3}
\end{table} 

(c) Let $C =\mathcal{G}(\BD_{2l(2k+1)}, \mathbb{Z}_{2k+1}, \BD_{2l(2p+1)}, \mathbb{Z}_{2p+1}, c_{g(w)} \circ \tau_{a,b})$ be such that $C \times D$ is not a subgroup of a freely acting splittable group. Then $l>2$ and  
$$C \times D= \mathcal{G}(\BD_{2l(2k+1)}, \mathbb{Z}_{2k+1}, \BD_{2l(2p+1)}, \mathbb{Z}_{2p+1}, \tau_{a,b}) \times \mathbb{Z}_{m'}$$
subject to the conditions 
$$4 \nmid m' \quad \wedge \quad \gcd (m',ms)=1.$$
Here the latter condition is required for both values of $$m=m_\varepsilon =\gcd (2p+1- \varepsilon (2k+1)a, 2l(2k+1), 2l(2p+1)(2k+1) ),$$ for  $\varepsilon = \pm 1$, and $$s= \frac{2l(2k+1)}{\gcd (2l(2k+1), 2l (2p+1)(2k+1))}.$$

(d) Let  $C=\mathcal{G}(\BD_{2(2k+1)}, \mathbb{Z}_{2k+1}, \mathbb{Z}_{4p}, \mathbb{Z}_p, \varphi(r))$ be such that $C \times D$ is not a subgroup of a freely acting splittable group. Then the latter belongs to the following table.
\begin{table}[H]
\begin{center}
\begin{tabular}{ ||c|c|| }
\hline
\multicolumn{2}{ ||c ||}{ $\mathcal{G}(\BD_{2(2k+1)}, \mathbb{Z}_{2k+1}, \mathbb{Z}_{4p}, \mathbb{Z}_p, \varphi(r))\times D, \quad \text{$p$ even}$ }  \\
\hline
$D$ &  Conditions \\ \hline
$\mathbb{Z}_m, \BD_{2m}$ & $\gcd (2k+1,p,m)=1$ \\
$\BT,\BO$ & $2,3 \nmid \gcd (2k+1,p)$ \\
$\BI$ & $2,3,5 \nmid \gcd (2k+1,p)$ \\ \hline
\end{tabular}
\end{center}
\end{table}

(e) Let $C= \mathcal{G}( \BD_{2(2k+1)}, \mathbb{Z}_{2k+1},  \BD_{2(2p+1)}  , \mathbb{Z}_{2p+1}, \varphi(r))$. Then $C \times D$ is a subgroup of a splittable group that acts freely on $M$.

(f) Let $C= \mathcal{G}(\mathbb{Z}_{2(2k+1)}, \mathbb{Z}_{2k+1}, \BO, \BT, \ID )$. Then either $3 \nmid (2k+1)$, in which case any group $C \times D$ acts freely on $M$, or in case $3|(2k+1)$, we have $D \in \{ \mathbb{Z}_m,\BD_{2m} \}$ subject to the condition $3 \nmid m$. 

(g)  Let $C= \mathcal{G}(\mathbb{Z}_{2(2k+1)}, \mathbb{Z}_{2k+1}, \BD_{4p} \BD_{2p}, \ID )$. Then $C \times D$ is a subgroup of a freely acting splittable group.

\et
\pf It suffices to figure out sufficient elements in 
$$\mathcal{W}(C)=\{\Real(x) \ : \ \Real(x)=\Real(y)  ,  (x,y)\in C \}\subset \mathbb{R}$$
 so we can read from them the precise conditions to impose on $C \times D $ to act freely on $M$.  
 
 (a) The maps defining $C = 
 \mathcal{G}(\mathbb{Z}_{kl}, \mathbb{Z}_{k}, \mathbb{Z}_{pl}  , \mathbb{Z}_p,  \varphi(r)) \subset \SU(2)^2$ are given by 
\begin{align}
\alpha & : \mathbb{Z}_{kl} \ni [x]_{kl}  \mapsto [rx]_l \in \mathbb{Z}_{l}\ , &   \beta : \mathbb{Z}_{pl} \ni [y]_{pl} \mapsto [y]_l \in \mathbb{Z}_{l}.\nonumber 
\end{align}
In particular, $C= \{ (e^{\frac{2 \pi i x}{kl} },e^{\frac{2 \pi i (rx+yl)}{pl} }) \ : \ y \in \mathbb{Z}\}$. To determine $\mathcal{W}(C)$ we must solve the  congruence system 
\begin{align}
px& = \varepsilon k(rx+yl) \mod lpk\nonumber 
\end{align}
separately for $\varepsilon \in \{-1,+1\}$. This leads to the following diophantine equation in three variables
$$ \frac{(p- \varepsilon kr)}{m}x- \frac{\varepsilon kl}{m}y+\frac{lpk}{m}z=0,$$
where $m=\gcd (p -\varepsilon  kr, kl,lpk)$.
Dividing this equation by $\frac{lk}{m}$ and setting $x=gx'$ for $x' \in \mathbb{Z} $, where $g=\gcd ( \frac{kl}{m}, \frac{lpk}{m})$, we obtain the equation
$$\frac{(p- \varepsilon kr)}{kl}gx'-\varepsilon y+pz=0,$$
which can be easily solved by considering  it as an inhomogeneous equation in the variables $y$ and $z$.  The general solution is given by
 \begin{align}
 x &=gx', \nonumber \\
 y &= \varepsilon \frac{p-\varepsilon kr}{kl} gx'+ \lambda p, \nonumber \\
 z &= \varepsilon \lambda, \nonumber 
 \end{align}
 where the parameters $x'$ and $\lambda$ are integers. Therefore  
 $$\mathcal{W}(C)=\Real (\mathbb{Z}_{ms}),\quad s =\frac{kl}{\gcd (kl,lpk)}.$$

(b) The maps defining the group $C = \mathcal{G}(\mathbb{Z}_{3n}, \mathbb{Z}_n,\BT,\BD_4,  \varphi(r)) \subset \SU(2)^2$ are given below.
\begin{align}
\alpha &: e^{\frac{2 \pi i x }{3n} } \mapsto e^{\frac{2 \pi i r x }{3} },  & 	\beta : z \mapsto \beta (z)= \left\{ 
  \begin{array}{l l}
    \text{$1$} & \quad \text{ if $z \in \BD_4$ }\\
    \text{$e^{\frac{2\pi i }{3} }$} & \quad \text{ if $z \in  \frac{-1+i+j+k}{2} \BD_4 $  } \\
		\text{$e^{\frac{4 \pi i  }{3} }$} & \quad \text{  if $z \in  \frac{1+i+j+k}{2} \BD_4. $}
  \end{array} \right . \nonumber 
	\end{align}
	Thus, for $r=1$ we get
\begin{align} 
	C= \biggl\{(e^{\frac{2 \pi i x}{n}},z_0),(e^{\frac{2\pi i(1+3y) }{3n} },z_1), (e^{\frac{2 \pi i(2+3z)  }{3n} },z_2) \ : \  z_0 \in \BD_4, \  z_1 \in  \frac{-1+i+j+k}{2} \BD_4, \nonumber 
\\ z_2 \in  \frac{1+i+j+k}{2} \BD_4, \  x,y,z \in \mathbb{Z} \biggl\} ,\nonumber
     \end{align}
     whereas for $r=2$, 
     \begin{align} 
	C= \biggl\{(e^{\frac{2 \pi i x}{n}},z_0),(e^{\frac{2\pi i(1+3y) }{3n} },z_1), (e^{\frac{2 \pi i(2+3z)  }{3n} },z_2) \ : \ z_0 \in \BD_4, \  z_1 \in  \frac{1+i+j+k}{2} \BD_4, \nonumber 
\\ z_2 \in  \frac{-1+i+j+k}{2} \BD_4, \  x,y,z \in \mathbb{Z} \biggl\} ,\nonumber
     \end{align}
     For any value of $r \in \{1,2\} $, $0$ or $-1$ is an element in $\mathcal{W}(C)$ precisely if $4 |n$ or $2|n$ respectively.
      Moreover, to verify if either $\frac{1}{2}\in \mathcal{W}(C)$ or $-\frac{1}{2}\in \mathcal{W}(C)$ amounts to solve each of the following equations, 
     \begin{align}
     \cos\left( \frac{2 \pi (1+3x)}{3n}\right) &= \pm \frac{1}{2} \ , \  \cos\left( \frac{2 \pi (2+3y)}{3n}\right)  = \pm \frac{1}{2}. \nonumber 
     \end{align}
     Our analysis must distinguish according to the parity of $n$. For $n$ even, the element $-1$ lies in $\mathcal{W}(C)$ for $r=1,2$, which implies  $D= \mathbb{Z}_{2p+1}$ in order to have a group acting freely. In such case $\frac{1}{2} \notin \mathcal{W}(C) \cap \Real(D)$, see Table \ref{tabla2}. Observe also that if $-\frac{1}{2} \in \mathcal{W}(C)$ and the group $C \times D $ acts freely on $M$, then similarly $3 \nmid (2p+1)$, and so  the former group will be a subgroup of the freely acting splittable group $A\times B \times D$. In the opposite case, i.e.\ $-\frac{1}{2} \not\in \mathcal{W}(C)$, we get a new freely acting semi-splittable group, see Lemma \ref{res}. 

Now, in case $n$ is odd, then $\mathcal{W}(C) \subset \{1, \pm \frac{1}{2}\}$. The precise elements in $\mathcal{W}(C)$ can be read off from Lemma \ref{res} and a case by case analysis leads to Table \ref{stt.3}.  

(c) Consider the group $C =\mathcal{G} ( \BD_{2l(2k+1)}, \mathbb{Z}_{2k+1}, \BD_{2l(2p+1)} , \mathbb{Z}_{2p+1}, c_{g(w)} \circ \tau_{a,b})$ for $l>2$, and its defining maps
\begin{align}
	 \alpha (z)= \left\{ 
  \begin{array}{l l}
    \text{$ c_{g(w)}(e^{\frac{i \pi x a}{l} }) $} & \quad \text{if $ z=e^{ \frac{i \pi x }{ l(2k+1) } } $ }\\
    \text{$ c_{g(w)}(je^{\frac{i \pi (a(x-1) +b+1)}{l} })$} & \quad \text{if $ z=je^{\frac{i \pi x }{ l(2k+1)}} $ } 	
  \end{array} \right. ,  \nonumber \\  \beta (z)= \left\{ 
  \begin{array}{l l}
    \text{$ e^{\frac{i \pi y }{l} } $} & \quad \text{if $ z=e^{ \frac{i \pi y }{ l(2p+1) } } $ }\\
    \text{$ je^{\frac{i \pi y}{l} }$} & \quad \text{if $ z=je^{\frac{i \pi y }{ l(2p+1)}}$ } 	
  \end{array} \right.  ,\nonumber
	\end{align}
where $(a,b) \in \mathbb{Z}_{2l}^{\times} \times \mathbb{Z}_{2l}$ and $g(w) \in \{e^{\frac{i \pi w}{l}},je^{\frac{i \pi w}{l} } \ : \ w \in \mathbb{Z} \}$. It is easy to see that zero is an element in $\mathcal{W}(C)$ for any $g(w) \in \BD_{2l}$, and so we see  the necessity of having $D= \mathbb{Z}_{m'}$ for an integer $m' \in \mathbb{Z}$ non-divisible by four. Furthermore, since $c_{g(w)}$ sends the element $e^{\frac{i \pi x a}{l} }$ either to itself or its inverse, the rest of the elements in $\mathcal{W}(C)$ different from zero are determined similarly as in (a). 
In fact, we have to solve
\[ Px= \varepsilon (Rx + Ly)(2k+1)\mod LP(2k+1),\]
where $P= 2p+1$, $L=2l$, $R=\varepsilon' a$ and $\varepsilon' =\varepsilon' (g(w)$ is as follows 
\begin{align}
	\varepsilon' (g(w))= \left\{ 
  \begin{array}{l l}
    \text{$1 $} & \quad \mbox{if} \quad \text{$ g(w)=e^{ \frac{i \pi x }{ l   }}$}\\
    \text{$ -1$} & \quad \mbox{if} \quad \text{$ g(w)=je^{\frac{i \pi x }{ l }}$} 	
  \end{array} \right. \nonumber
	\end{align}  
Summarizing, we obtain 
$$\mathcal{W}(C)= \{0\} \cup \Real(\mathbb{Z}_{m s}) ,$$
where $$m=\gcd (2p+1- \varepsilon \varepsilon' (2k+1)a, 2l(2k+1), 2l(2p+1)(2k+1) ), \quad \varepsilon \in \{ -1, 1\},$$ and $$s= \frac{2l(2k+1)}{\gcd (2l(2k+1), 2l (2p+1)(2k+1))}.  $$
This leads to the conditions stated in the theorem after absorbing the sign $\epsilon'$ into $\epsilon$. 

 Now, let $l=2$ and recall that any automorphism of $\BD_4$ is obtained by conjugation $\varphi_w$  with an element $w\in \BO$. The maps defining $C $ are:  
 \begin{align}
	\alpha (z)& = \left\{ 
  \begin{array}{l l}
  \text{$ 1 $} & \quad \text{if $ z=e^{ \frac{2 \pi ix }{ (2k+1) } } $ }\\
    \text{$ \varphi_w(j)$} & \quad \text{if $ z=je^{\frac{2 \pi ix }{ (2k+1)}} $ }  \\
    \text{$ \varphi_w(i) $} & \quad \text{if $ z=e^{ \frac{i \pi (1+4x) }{ 2(2k+1) } } $ }\\
    \text{$ -\varphi_w(k)$} & \quad \text{if $ z=je^{\frac{i \pi (1+4x) }{ 2(2k+1)}} $ }  \\
    \text{$ -1 $} & \quad \text{if $ z=e^{ \frac{i \pi (1+2x) }{ (2k+1) } } $ }\\
    \text{$ -\varphi_w(j)$} & \quad \text{if $ z=je^{\frac{i \pi (1+2x) }{ (2k+1)}} $ } \\
    \text{$ -\varphi_w(i) $} & \quad \text{if $ z=e^{ \frac{i \pi (3+4x) }{ 2(2k+1) } } $ }\\
    \text{$ \varphi_w(k)$} & \quad \text{if $ z=je^{\frac{i \pi (3+4x) }{ 2(2k+1)}} $ } 	
    \end{array} \right. \nonumber &
  \beta (z)& = \left\{ 
  \begin{array}{l l}
  \text{$ 1 $} & \quad \text{if $ z=e^{ \frac{2 \pi iy }{ (2p+1) } } $ }\\
    \text{$ j$} & \quad \text{if $ z=je^{\frac{2 \pi iy }{ (2p+1)}} $ }  \\
    \text{$ i $} & \quad \text{if $ z=e^{ \frac{i \pi (1+4y) }{ 2(2p+1) } } $ }\\
    \text{$ -k$} & \quad \text{if $ z=je^{\frac{i \pi (1+4y) }{ 2(2p+1)}} $ }  \\
    \text{$ -1 $} & \quad \text{if $ z=e^{ \frac{i \pi (1+2y) }{ (2p+1) } } $ }\\
    \text{$ -j$} & \quad \text{if $ z=je^{\frac{i \pi (1+2y) }{ (2p+1)}} $ } \\
    \text{$ -i $} & \quad \text{if $ z=e^{ \frac{i \pi (3+4y) }{ 2(2p+1) } } $ }\\
    \text{$ k$} & \quad \text{if $ z=je^{\frac{i \pi (3+4y) }{ 2(2p+1)}} $ } 	
   \end{array} \right. \nonumber
	\end{align}    
Observe that for any choice of $w$ we have that $\alpha(e^{i \pi })=\beta( e^{i \pi })=-1$, and so $D= \mathbb{Z}_{m}$ with $2 \nmid m$. Moreover, since 
$$\alpha (e^{ \frac{2 \pi ix }{ (2k+1) } })= \beta (e^{ \frac{2 \pi iy }{ (2p+1) } }) =1,$$
we have $\mathbb{Z}_{2k+1} \times \mathbb{Z}_{2p+1} \subset C$.
 It follows that $\Real(\mathbb{Z}_{\gcd (2k+1,2p+1)} ) \subset \mathcal{W}(C)$, and hence the necessity to have  
$$1=\gcd (2k+1,2p+1,m)=\gcd (4(2k+1),4(2p+1),m).$$

(d) Now, consider $C= \mathcal{G}(\BD_{2(2k+1)}, \mathbb{Z}_{2k+1},  \mathbb{Z}_{4p}, \mathbb{Z}_p, \varphi(r))$. The group $C $ is equal to one of the following groups:
\begin{align*}
 \biggl\{ (e^{\frac{2\pi i x}{2k+1} },e^{\frac{2\pi i y}{p} }),(je^{\frac{2 \pi i x}{2k+1} },e^{\frac{ \pi i (1+4y)}{2p} }),(e^{\frac{i \pi x}{2k+1} }, e^{\frac{i \pi(1+2y)}{p} }),  (je^{\frac{i \pi x}{2k+1} }   , e^{\frac{i \pi (3+4y)}{2p} })    \ : \  x,y \in \mathbb{Z}  \biggl\}, \\
 \biggl\{ (e^{\frac{2\pi i x}{2k+1} },e^{\frac{2\pi i y}{p} }),(je^{\frac{i \pi x}{2k+1} },e^{\frac{ \pi i (1+4y)}{2p} }),(e^{\frac{i \pi x}{2k+1} }, e^{\frac{i \pi(1+2y)}{p} }),  (je^{\frac{2 \pi i x}{2k+1} }   , e^{\frac{i \pi (3+4y)}{2p} })    \ : \  x,y \in \mathbb{Z}  \biggl\}.
\end{align*} 
Let us distinguish according to the parity of $p$. If $p$ is odd, then for any choice of $\varphi(r)$ we get $-1 \in \mathcal{W}(C)$ and also 
$$\Real(\mathbb{Z}_{\gcd (2k+1,p)})\subset \mathcal{W}(C).$$ 
In order for $C \times D$ to act freely on $M$, it is necessary to have $D= \mathbb{Z}_m$ with $2 \nmid m$ and
$$1=\gcd (2k+1,p,m)= \gcd (2(2k+1),4p,m).$$
This shows that $C \times D $ is a subgroup of a freely acting splittable group. Let now $p$ be even. In this case $0 \notin \mathcal{W}(C)$ and $$\mathcal{W}(C) \supset \Real (\mathbb{Z}_{\gcd (2k+1,p)}).$$ 
The rest of the elements in $\mathcal{W}(C)$ are obtained by solving the following equations
$$px=\varepsilon (2k+1)(1+2y) \mod 2p(2k+1),$$
for any $\varepsilon \in \{1,-1\}$. Because $p$ is even, these equations have no solution. In other words, $\mathcal{W}(C) =\Real (\mathbb{Z}_{\gcd (2k+1,p)})$ . 

(e) The group $C= \mathcal{G}(\BD_{2(2k+1)}, \mathbb{Z}_{2k+1}, \BD_{2(2p+1)}, \mathbb{Z}_{2p+1}, \varphi(r))$ is either one of the following groups
\begin{align*}
 \biggl\{ (e^{\frac{2\pi i x}{2k+1} },e^{\frac{2\pi i y}{2p+1} }),j(e^{\frac{2 \pi i x}{2k+1} },e^{\frac{ 2\pi i y }{2p+1} }),(e^{\frac{i \pi (1+2x)}{2k+1} }, e^{\frac{i \pi(1+2y)}{2p+1} }),  j(e^{\frac{i \pi (1+2x)}{2k+1} }   , e^{\frac{i \pi (1+2y)}{2p+1} })    : \  x,y \in \mathbb{Z}  \biggl\}, \\
 \biggl\{ (e^{\frac{2\pi i x}{2k+1} },e^{\frac{2\pi i y}{2p+1} }),j(  e^{\frac{i \pi (1+2x)}{2k+1} }        ,e^{\frac{ 2\pi i y }{2p+1} }),(e^{\frac{i \pi (1+2x)}{2k+1} }, e^{\frac{i \pi(1+2y)}{2p+1} }),  j(e^{\frac{2 \pi i x}{2k+1} }   , e^{\frac{i \pi (1+2y)}{2p+1} })     : \  x,y \in \mathbb{Z}   \biggl\}.
\end{align*} 
As in the first part of (d), we see that $-1 \in \mathcal{W}(C)$ and $\mathcal{W}(C) \supset \mathbb{Z}_{\gcd (2k+1,2p+1)}$. Therefore, we must have $D=\mathbb{Z}_m$ with $m$ odd and $$1=\gcd (2k+1,2p+1,m)=\gcd (2(2k+1),2(2p+1),m).$$
Again this shows that $C \times D$ is a subgroup of a freely acting splittable group.

(f) The group $C= \mathcal{G}(\mathbb{Z}_{2(2k+1)}, \mathbb{Z}_{2k+1}, \BO, \BT, \ID )$ is given by 
$$\left\{ \left(e^{\frac{2 \pi i x}{2k+1}}, z_1) \right),\left(e^{\frac{ \pi i (1+2x)}{2k+1}}, z_2\right) \ : \ z_1 \in \BT , z_2 \in e^{\frac{i \pi }{4}} \BT \right\}.$$
From this, one can verify that $\mathcal{W}(C) \subset \{ 1, -\frac{1}{2}\}$. This shows the claim.

(g) Let $C= \mathcal{G}(\mathbb{Z}_{2(2k+1)}, \mathbb{Z}_{2k+1}, \BD_{4p} \BD_{2p}, \ID )$. It is readily seen that $-1 \in \mathcal{W}(C)$ and that the latter set contains  $\Real (\mathbb{Z}_{\gcd (2(2k+1),2p)})$. We conclude the necessity to impose the conditions: $D= \mathbb{Z}_m$ with $m$ odd and   
$$1= \gcd ( 2( 2k+1), 2p,m)=  \gcd ( 2( 2k+1), 4p,m). $$
No new freely acting group is obtained in this way.
 \epf 
 The following theorem ends up our classification of freely acting finite subgroups $C \times D \subset \SU(2)^3$.
 \begin{Th}\label{q.final2} Let $C \times D $ be a type III freely acting group different from any group occurring in Theorem \ref{q.final}. Then, either $C \times D $  is a subgroup of a freely acting splittable group, or it is one of the following groups. 
 \begin{table}[H]
\begin{center}
\begin{tabular}{ ||c |c |c |c | c|| }
\hline
\multicolumn{5}{ ||c|| }{$ \mathcal{G}(A,A_0,B,B_0, \theta) \times \mathbb{Z}_m , \quad 2\nmid m $ }   \\ \hline
 $A$ & $A_0$ & $B$ & $B_0$ &Conditions  \\ \hline 
 $\BI$ & $\mathbb{Z}_2$ & $\BI$ & $\mathbb{Z}_2$ & $3 \nmid m$ \\ \hline
  $\BD_{6k}$& $\mathbb{Z}_{2k}$ & $\BO$ &$\BD_4$ &  $3 |k $ \\  \hline
  $\BD_{2kl}$& $\mathbb{Z}_{2k}$ & $\BD_{2pl}$& $\mathbb{Z}_{2p}$ & $\gcd (\tilde{m}s, m)=1$ \\  \hline 
\end{tabular}
\end{center}
\end{table}
Where the above conditions for the group $ \mathcal{G}(\BD_{2kl}, \mathbb{Z}_{2k},  \BD_{2pl}, \mathbb{Z}_{2p}, c_w \circ \tau_{a,b}) \times \mathbb{Z}_m$, are required for both values of $$\tilde{m}=\tilde{m}_\varepsilon=\gcd (p -\varepsilon  ka, kl,lpk)),$$ for $\varepsilon = \pm 1$, and $$ s =\frac{kl}{\gcd (kl,lpk)}.$$
\end{Th}
 \pf  We are left with the case $D = \mathbb{Z}_{m}$ for $m\in \mathbb{Z}$ odd, see Proposition \ref{form.strict}. It suffices again to find suitable elements in $\mathcal{W}(C)$. Since we have numerous choices for $C= \mathcal{G}(A,A_0,B,B_0,\theta) $, let us proceed by considering all possible candidates for $F=\sfrac{B}{B_0}$ according to Tables \ref{tabla-p} and \ref{tabla-p.1}, which were not considered in Theorem \ref{q.final}. 

(a) Let $F \in \{\T,\OR,\I\}$. In those cases $A=B \in \{\BT,\BO, \BI\}$ and $A_0=B_0=\mathbb{Z}_2$. Let $\alpha: A \mapsto F$ and $\beta: B \mapsto F$ be the maps defining $C$. Since the automorphisms $\theta$ of $\sfrac{B}{\mathbb{Z}_2}$, for $B = \BT, \BO$, is induced by conjugation with a class $[w]\in\sfrac{\BO}{\mathbb{Z}_2}$ , we have  $\alpha(x)=\beta(c_w(x))$ for any element $x \in B$. It follows that $\mathcal{W}(C)= \Real(B)$. This shows that $C \times D$ is a subgroup of a freely acting splittable group.

We treat now the case $B= \BI$. Because of the argument above, it suffices to consider the (only) non-trivial outer automorphism of $F$ in the construction of $C$. Since two elements $[z],[w]$ in $\sfrac{\BI}{\mathbb{Z}_2}$ are conjugate (precisely) when $\Real(z)=\pm \Real(w)$, see Table \ref{tabla6}, the group $\sfrac{\BI}{\mathbb{Z}_2}$ has five conjugacy classes $\mathcal{C}([y_i])$, $i=1,...,5$, with representatives $y_i \in \BI$ having (up to sign) every possible real part of an element in $\BI$. On the other hand, if $\varphi: \sfrac{\BI}{\mathbb{Z}_2}\mapsto \sfrac{\BI}{\mathbb{Z}_2}$ is the representative of the generator of $\Out(\sfrac{\BI}{\mathbb{Z}_2})$ described in item (6)  Section \ref{automorphisms}, then it fixes all conjugacy classes except of the two conjugacy classes of order 12, which are exchanged. It follows that $\mathcal{W}(C)=\{0 , \pm 1, \pm \frac{1}{2}\}$. This leads to the additional condition $3 \nmid m$.

(b) Groups $\mathcal{G}(A,A_0,B,B_0, \theta)  \subset \SU(2)^2$ giving rise to $F \cong \BD_{2l}$ were analyzed in Theorem \ref{q.final}.  

(c) Let $F\cong \D_{2l}$ for some integer $l\geq 2$. Let $C=\mathcal{G}(\BO, \BD_4,\BO, \BD_4, c_{w} \circ \tau_{a,b} )$, where $(a,b)\in \mathbb{Z}_6^{\times} \times \mathbb{Z}_6$ and $w\in \BO $. We easily check that 
$$\left(- \frac{1}{2}(1+i+j+k), c_{w}\left(-\frac{1}{2}(1+i+j+k) \right)  \right)\in C,$$ 
whence $a  = 1 \mod 6$. As for if $a=5 \mod 6$, then 
$$\left(- \frac{1}{2}(1+i+j+k), c_{ w }\left(\frac{1}{2}(-1+i+j+k) \right)  \right)\in C.$$
 This implies that $- \frac{1}{2} \in \mathcal{W}(C)$ for any homomorphism $c_w \circ \tau_{a,b}$. This shows that $C \times D$ is a subgroup of a freely acting splittable group.

 Now, consider $C=\mathcal{G}(\BD_{6p}, \mathbb{Z}_{2p}, \BO, \BD_4  , c_{g(w)} \circ \tau_{a,b} )$, where $(a,b) \in \mathbb{Z}_6^{\times} \times \mathbb{Z}_6$ and $g(w)\in  \BD_{6p} =\{e^{\frac{i \pi w}{3p} }, je^{\frac{i \pi w}{3p} }: w \in \mathbb{Z}\}$. For $a=1 \mod 6$, we have 
 \begin{align}
	 \alpha \left( e^{\frac{2 \pi i}{3} }\right)= \left\{ 
  \begin{array}{l l}
    \text{$  \left( \frac{1}{2}(1+i+j+k) \right)^{2p} \BD_4$} & \quad \text{$ g(w)= e^{\frac{i \pi w}{3p} }$ }\\
    \text{$  \left(\frac{1}{2}(1+i+j+k) \right)^{4p} \BD_4$} & \quad \text{$ g(w)=je^{\frac{i \pi w}{3p} } $  ,} 
  \end{array} \right.  \nonumber
	\end{align}
 whereas for $a=5 \mod 6$   
 \begin{align}
	 \alpha \left( e^{\frac{2 \pi i}{3} }\right)= \left\{ 
  \begin{array}{l l}
    \text{$  \left( \frac{1}{2}(1+i+j+k)  \right)^{4p} \BD_4$} & \quad \text{$ g(w)= e^{\frac{i \pi w}{3p} }$ }\\
    \text{$ \left(  \frac{1}{2}(1+i+j+k)  \right)^{2p} \BD_4$} & \quad \text{$ g(w)=je^{\frac{i \pi w}{3p} } $ , } 
  \end{array} \right.   \nonumber
	\end{align}
	where $\beta: \BO \mapsto \sfrac{\BO}{\BD_4}$ is in both cases the natural map. Using the fact that the order of the element
	$$ \frac{1}{2}(1+i+j+k) \BD_4 \in \sfrac{\BO}{\BD_4} $$
	is $3$, it follows that $- \frac{1}{2} \in \mathcal{W}(C)$ precisely if $3 \nmid p$. In that case $C \times D$ is a subgroup of a freely acting splittable group. Otherwise $\mathcal{W}(C) \cap \Real( D)=\{1\}$, and so the group $$\mathcal{G}(\BD_{6p}, \mathbb{Z}_{2p}, \BO,  \BD_4  , c_{g(w)} \circ \tau_{a,b} ) \times D $$ acts freely on $M$.
	
	Lastly, consider the case in which $C= \mathcal{G}(\BD_{2kl},\mathbb{Z}_{2k},\BD_{2pl},\mathbb{Z}_{2p}, c_w \circ \tau_{a,b})$.
	To determine the elements of $\mathcal{W}(C)$ different from zero, it suffices to consider the subset of $C \subset \SU(2)^2$ given by 
				$$\{(e^{\frac{i \pi x}{kl}}, e^{\frac{i \pi ( ax+l  y)}{pl}}) \ : \ x,y \in \mathbb{Z} \},$$
		and so solve the equation  
		$$px=k \varepsilon (ax+l  y) \mod 2kpl. $$
		In fact, because $m$ is odd, and hence $$\gcd (2kl,2pl,m)= \gcd (kl,pl,m),$$ we can consider 
		$$px=k \varepsilon (ax+l  y) \mod kpl $$
		instead. This situation was encountered in part (a) of Theorem \ref{q.final}, which applied to this situation shows the claim.

(d) Let $F \cong \mathbb{Z}_{2}$. For the group $ C=\mathcal{G}(\BO,\BT,\BO,\BT , \ID)$, we have $\mathcal{W}(C) = \Real(\BO)$, and so no new freely acting group will be obtained. 

Consider the group $C= \mathcal{G}(\BO, \BT, \mathbb{Z}_{2k}, \mathbb{Z}_k, \ID)$ and let
\begin{align}
	\alpha &: z \mapsto \alpha (z)= \left\{ 
  \begin{array}{l l}
    \text{$1$} & \quad \text{if $z \in \BT$ }\\
    \text{$-1$} & \quad \text{if $z \in e^{\frac{i \pi}{4} }\BT $ } 	
  \end{array} \right. , \  \beta \left( e^{\frac{i \pi x}{k} }\right) =  \left\{ 
  \begin{array}{l l}
    \text{$1$} & \quad \text{if $ x = 0 \mod 2$ }\\
    \text{$-1$} & \quad \text{if $ x = 1 \mod 2 $ } \end{array} \right.  \nonumber
	\end{align}
	 be the maps defining it. We check that $-\frac{1}{2}$ lies in $\mathcal{W}(C)$ precisely if $3|k$. It follows the necessity to impose on $C \times D$ the condition $3 \nmid  \gcd (k,m)$, which implies that $3 \nmid  \gcd (2k,m)$, and so no new freely acting group will be obtained in this fashion.
	 
	  The group $C= \mathcal{G}( \BD_{4k}, \BD_{2k}, \BO, \BT, \ID)$ is given by
\begin{align*}
\biggl\{ (e^{\frac{i \pi x}{k} }, z),(je^{\frac{i \pi x}{k} },z),(e^{\frac{i \pi (2x+1)}{2k} }, e^{\frac{i \pi}{4} }z),  (je^{\frac{i \pi (2x+1)}{2k} }   , e^{\frac{i \pi}{4} }z)   \ : \ z\in \BT, x\in \mathbb{Z}    \biggl\}.
\end{align*} 
A necessary condition to impose on $C \times D $ to act freely on $M$ reads $ 3 \nmid \gcd (k,m)= \gcd (4k,m)$. Therefore, no new freely acting subgroup will be found in this way.

 The group $C= \mathcal{G}(\mathbb{Z}_{2k}, \mathbb{Z}_k, \BD_{4p}, \BD_{2p}, \ID)$ is given as follows 
\begin{align*}
\biggl\{ (e^{\frac{2i \pi x}{k} },e^{\frac{i \pi y}{p} }),(e^{\frac{2i \pi x}{k} },je^{\frac{i \pi y}{p} }),(e^{\frac{i \pi (2x+1)}{k} },e^{\frac{i \pi (2y+1)}{2p} } ), (e^{\frac{i \pi (2x+1)}{k} } 
  ,je^{\frac{i \pi (2y+1)}{2p} } ) \ : \ x, y \in \mathbb{Z}  \biggl\}.
\end{align*} 
In consequence, $\mathcal{W}(C)= \Real(\mathbb{Z}_{\gcd (k,2p)})$. Therefore the necessity to impose 
$$1= \gcd (k,2p,m)= \gcd (2k,4p,m),$$
and so there is no new freely acting subgroup.

 Analogously, for any group $C=\mathcal{G}(\BD_{4k}, \BD_{2k}, \BD_{4p},  \BD_{2p}, \ID)$, we have to have  $$1= \gcd (4k,4p,m),$$ and so no new freely acting subgroup is obtained in this fashion.

(e) Let $F \cong \mathbb{Z}_3$. The analysis for $\mathcal{G}(\mathbb{Z}_{3k},\mathbb{Z}_{k},\BT ,\BD_4,\varphi(r))$ was performed in part (b) of Theorem \ref{q.final}, whereas for $C=\mathcal{G}(\BT, \BD_4, \BT,\BD_4, \varphi(r))$, it is not difficult to see that $\mathcal{W}(C)=\Real(\BT)$, and so no new freely acting groups are obtained in this case.	
	
	(f) The groups $C \subset \SU(2)^2$ leading to $F \cong \mathbb{Z}_l$ for $l \notin \{2,3\}$ were analyzed in Theorem \ref{q.final}.
	
\epf

 
\end{document}